\documentstyle[12pt,verbatim]{amsart}
\let\widehat\hat
\let\mathcal\cal
\newtheorem{theorem}{Theorem}[section]

\newtheorem{lemma}[theorem]{Lemma}
\newtheorem{corollary}[theorem]{Corollary}
\newtheorem{fact}[theorem]{Fact}

\newtheorem{case}{Case}

\theoremstyle{definition}

\newtheorem{definition}[theorem]{Definition}

\theoremstyle{remark}
\newtheorem{remark}{Remark}

\newcommand{\proof}{\begin{pf}}
\newcommand{\Proof}[1]{\begin{pf*}{Proof of #1}}
\newcommand{\eproof}{\end{pf}}
\newcommand{\Eproof}{\end{pf*}}
\newcommand{\sproof}[1]{\begin{pf*}{#1}}
\newcommand{\esproof}{\end{pf*}}

\newcommand{\arablabel}{
          \renewcommand{\labelenumi}{{\rm (\arabic{enumi})}}
          \renewcommand{\theenumi}{{\rm (\arabic{enumi})}}
          \renewcommand{\labelenumii}{{\rm (\arabic{enumii})}}
          \renewcommand{\theenumii}{{\rm (\arabic{enumii})}}
                    }

\newcommand{\alabel}{
          \renewcommand{\labelenumi}{{\rm (\alph{enumi})}}
          \renewcommand{\theenumi}{{\rm (\alph{enumi})}}
          \renewcommand{\labelenumii}{{\rm (\alph{enumii})}}
          \renewcommand{\theenumii}{{\rm (\alph{enumii})}}
                    }
\newcommand{\Alabel}{
          \renewcommand{\labelenumi}{{\rm (\Alph{enumi})}}
          \renewcommand{\theenumi}{{\rm (\Alph{enumi})}}
          \renewcommand{\labelenumii}{{\rm (\Alph{enumii})}}
          \renewcommand{\theenumii}{{\rm (\Alph{enumii})}}
                    }
\newcommand{\rlabel}{
          \renewcommand{\labelenumi}{{\rm (\roman{enumi})}}
          \renewcommand{\theenumi}{{\rm (\roman{enumi})}}
          \renewcommand{\labelenumii}{{\rm (\roman{enumii})}}
          \renewcommand{\theenumii}{{\rm (\roman{enumii})}}
                    }

\def\myheads#1;#2;{
\pagestyle{myheadings}
\markboth{{\sc\hfill #1\hfill\protect\makebox[0cm][r]{\rm\today}}}
{{\sc\protect\makebox[0cm][l]{\rm\today}\hfill #2\hfill}}
}

\newcommand{\acal}{{\mathcal A}}
\newcommand{\bcal}{{\mathcal B}}
\newcommand{\ccal}{{\mathcal C}}
\newcommand{\dcal}{{\mathcal D}}

\newcommand{\fcal}{{\mathcal F}}

\newcommand{\ical}{{\mathcal I}}

\newcommand{\mcal}{{\mathcal M}}
\newcommand{\ncal}{{\mathcal N}}
\newcommand{\pcal}{{\mathcal P}}

\newcommand{\scal}{{\mathcal S}}
\newcommand{\tcal}{{\mathcal T}}
\newcommand{\ucal}{{\mathcal U}}
\newcommand{\vcal}{{\mathcal V}}

\newcommand{\xcal}{{\mathcal X}}

\newcommand{\setm}{\setminus}
\newcommand{\empt}{\emptyset}
\newcommand{\subs}{\subset}
\newcommand{\sups}{\supset}
\newcommand{\oo}{{{\omega}_1}}
\newcommand{\rest}{\lceil}
\newcommand{\dom}{\operatorname{dom}}
\newcommand{\ran}{\operatorname{ran}}
\def\<{\left\langle}
\def\>{\right\rangle}
\def\cf{\operatorname{cf}}

\def\OO{{\omega}}
\def\oo{\omega_1}
\def\br#1;#2;{\bigl[ {#1} \bigr]^ {#2} }
\def\bc#1;#2;{\bigl( {#1} \bigr)^ {#2} }

\def\ooseq#1;#2;{\< {#1}_{#2}:{#2}<\oo\>}
\def\ooset#1;#2;{\{ {#1}_{#2}:{#2}<\oo\}}
\def\seq#1;#2;#3;{\< {#1}_{#2}:{#2}<#3\>}
\def\set#1;#2;#3;{\{ {#1}_{#2}:{#2}<#3\}}
\def\oseq#1;#2;{\< {#1}_{#2}:{#2}<\OO\>}
\def\oset#1;#2;{\{ {#1}_{#2}:{#2}<\OO\}}
\def\oosequ#1;#2;{\< {#1}^{#2}:{#2}<\oo\>}
\def\oosetu#1;#2;{\{ {#1}^{#2}:{#2}<\oo\}}
\def\sequ#1;#2;#3;{\< {#1}^{#2}:{#2}<#3\>}
\def\setu#1;#2;#3;{\{ {#1}^{#2}:{#2}<#3\}}
\def\osequ#1;#2;{\< {#1}^{#2}:{#2}<\OO\>}
\def\osetu#1;#2;{\{ {#1}^{#2}:{#2}<\OO\}}

\def\force{\raisebox{1.5pt}{\mbox {$\scriptscriptstyle\|$}}
\mbox{$\!\mbox{---}$}}

\newcommand{\notforce}{\not\!\!\!\force}

\newcommand{\fn}{\operatorname{Fn}}

\def\to{\longrightarrow}

\newcommand{\newcases}{\setcounter{case}{0}}

\def\fin#1;{\br #1;<{\omega};}

\newcommand\adot{{\dot A}}
\newcommand\bdot{{\dot B}}

\newcommand\dotd{{\dot D}}

\newcommand{\w}{\operatorname{w}}

\def\oot{{{\omega}_2}}
\def\con#1;{\mbox{\rm Con}$($#1$)$}

\def\qstar#1;{q^*_{#1}}
\def\nea{neighbourhood assignment }

\newcommand{\shat}{\operatorname{\hat s}}
\newcommand{\hhat}{\operatorname{\hat h}}
\newcommand{\zhat}{\operatorname{\hat z}}

\def\diamin#1;{D({#1})}
\def\dia#1;{D^s({#1})}
\newcommand{\fab}{{\varphi}_{{\alpha},{\beta}}}

\def\aaa#1;#2;{A({#2},#1)}
\def\aaadot#1;#2;{\dot{A}({#2},#1)}
\def\bbb#1;#2;{B(#2,#1)}
\def\bbbdot#1;#2;{\dot{B}({#2},{#1})}
\def\ddd#1;#2;{D(#2,#1)}
\def\coh#1;{\ccal_{#1}}
\def\lsups#1;#2;{\widehat{#1}^{\it s}}
\def\lsup#1;#2;{\widehat{#1}}
\def\ccc#1;{C^s({#1})}
\def\ccci#1;{\widehat{C}^s({#1})}
\def\cccmin#1;{C({#1})}
\def\cccimin#1;{\widehat{C}({#1})}
\def\cccc#1;#2;{C^{(#2)}({#1})}
\def\eee#1;{F({#1})}
\def\eeemin#1;{F^-({#1})}
\def\eee#1;{F^s({#1})}
\def\eeemin#1;{F({#1})}

\def\pii#1;#2;#3;{\pi^{#1}}

\newcommand{\supp}{\operatorname{supp}}
\newcommand{\sprd}{\operatorname{s}}
\newcommand{\hl}{\operatorname{h}}

\def\ii#1;#2;{i(#1,#2)}

\newcommand{\nok}{\<n_0,\dots,n_{k-1}\>}

\newcommand{\htt}{\operatorname{ht}}

\def\injs#1;#2;{(#1)^{#2}}
\def\injc#1;{(#1)}
\def\aint#1;#2;{A(#1,#2)}
\def\bint#1;#2;{B(#1,#2)}
\def\bcaldot{{\dot\bcal}}

\def\ig#1;#2;{\operatorname{val}(#2,#1)}
\def\matr#1;{\mcal(#1)}
\def\matrn#1;{\ncal(#1)}

\def\pint#1;#2;{#1_{(#2)}}

\def\lev#1;#2;{#1^{(#2)}}
\def\levup#1;#2;{#1^{[#2]}}
\def\levx#1;{\lev X;#1;}
\def\levxup#1;{\levup X;#1;}

\theoremstyle{plain}
\newtheorem{ertheorem}{Erd{\H os}-Rado Theorem}

\author{I. Juh\'asz}
\address{Mathematical Institute of the Hungarian Academy of Sciences}
\email{juhasz@@math-inst.hu}
\author{L. Soukup}
\address{Mathematical Institute of the Hungarian Academy of Sciences}
\thanks{The second author  was supported by DFG (grant Ko 490/7-1)}
\email{soukup@@math-inst.hu}
\author{Z. Szentmikl{\'o}ssy}
\address{E{\"o}tv{\"o}s Lor{\'a}nd University, Department of Analysis}
\email{szetmiklossy@@math-inst.hu}
\title{Combinatorial principles from adding Cohen reals}
\thanks{The preparation of this paper was supported by the 
Hungarian National Foundation for Scientific Research grant no. 16391.}

\begin{document}
\maketitle

\begin{abstract}
We first formulate several ``combinatorial principles''
concerning ${\kappa}\times {\omega}$ matrices of subsets of
${\omega}$ and prove that they are valid in the generic extension 
obtained by adding any number of
Cohen reals to any ground model $V$, provided that the parameter 
${\kappa}$ is an
${\omega}$-inaccessible regular cardinal in $V$.

Then in section \ref{sc:applications}
we present a large number of applications of these principles, 
mainly to topology.  Some of these consequences had been  established earlier in generic extensions obtained by adding ${\omega}_2$ Cohen reals 
to ground models satisfying $CH$, mostly  for the case ${\kappa}={\omega}_2$.
\end{abstract}


\myheads Combinatorial principles $\dots$;Combinatorial principles $\dots$;

\section{Introduction}\label{sc:int}
The last 25 years have seen a furious activity in proving results that are 
independent of the usual axioms of set theory, that is ZFC.
As the methods of these independence proofs 
(e.g. forcing or the fine structure theory of the constructible universe) 
are often rather sophisticated, while the results themselves are usually
of interest to ``ordinary'' mathematicians 
(e.g. topologists or analysts), it has been natural to try to isolate a relatively small number of principles, i.e.
independent statements that a) are {\em simple} to formulate and 
b) are {\em useful} in the sense that they have many interesting consequences.
Most of these statements, we think by necessity, are of combinatorial nature,
hence they have been called combinatorial principles.

In this paper we propose to present several new
 combinatorial principles that are all statements about $\pcal({\omega})$, the 
power set of the natural numbers. In fact, they all concern
matrices of the form
$\<\aaa n;{\alpha};:\<{\alpha},n\>\in {\kappa}\times {\omega}\>$, where
$\aaa n;{\alpha};\subs {\omega}$ for each 
$\<{\alpha},n\>\in {\kappa}\times {\omega}$, and, in the interesting cases,
${\kappa}$ is a regular cardinal with $c=2^{\omega}\ge{\kappa}>\oo$.

We show that these statements are valid in the generic extensions 
obtained by adding any number of Cohen reals to any ground model $V$,
assuming that the parameter ${\kappa}$ is a regular and 
${\omega}$-inaccessible cardinal in $V$ ( i.e.
${\lambda}<{\kappa}$ implies ${\lambda}^{\omega}<{\kappa}$).

Then we present a large number of consequences of these principles, 
some of them combinatorial but most of them topological,
mainly concerning separable and/or countably tight topological spaces.
(This, of course, is not  surprising because these are objects whose 
structure depends basically on $\pcal({\omega})$.)

The above  formulated criteria a) and b) as to what constitutes a
combinatorial principle are often contrary to each other:
for more usefulness one often has to sacrifice some simplicity.
It is not clear whether an ideal balance exists between them.
It is up to the reader to judge if we have come close to this balance.

\section{The combinatorial principles}\label{sc:principles}

The principles we formulate here are all statements on 
${\kappa}\times {\omega}$ matrices of subsets of
${\omega}$ claiming -- roughly speaking -- that all these matrices contain
large ``submatrices'' satisfying certain homogeneity properties.

To simplify the formulation of our results we introduce
the following pieces of notation. 
If $S$ is an arbitrary set and 	$k$ is a natural number
then let 	
$$
\injs S;k;=\{s\in S^k:|\ran s|=k\}
$$
and 
$$
\injs S;<{\omega};=\bigcup_{k<{\omega}}\injs S;k;.
$$
For $D_0,\dots, D_{k-1}\subs S$ we let
$$
\injc D_0,\dots,D_{k-1};=\{s\in\injs S;k;:\ \forall i\in k \ (s(i)\in D_i)\}.
$$

\begin{definition}\label{df:matrix}
If $S$ is a set of ordinals denote by 
$\matr S;$ the family of all $S\times{\omega}$-matrices
of subsets of ${\omega}$, that is, $\acal\in\matr S;$
if and only if $\acal=\<\aaa i;{\alpha};:{\alpha}\in S,i<{\omega}\>$,
where $\aaa i;{\alpha};\subs {\omega}$ for each ${\alpha}\in S$ 
and $i<{\omega}$. If 
$\acal=\<\aaa i;{\alpha};:{\alpha}\in S,i<{\omega}\>\in\matr S;$
and $R\subs S$ we define the restriction  of $\acal$ to $R$, $\acal\rest R$ 
in the straightforward way: 
$\acal\rest R=\<\aaa i;{\alpha};:{\alpha}\in R,i<{\omega}\>$.
If 
$\acal=\<\aaa i;{\alpha};:{\alpha}\in S,i<{\omega}\>\in \matr S;$,
$t\in\omega^{<{\omega}}$ and $s\in \injs S;|t|;$   then we let
$$
\aint s;t;=\bigcap_{i<|t|}\aaa t(i);s(i);.
$$
\end{definition}

Now we formulate our first and probably most important principle
that we call $\ccc {\kappa};$. 
We also specify a weaker version of $\ccc {\kappa};$ denoted
by $\cccmin {\kappa};$ because in most of the applications 
(\ref{tm:tower}, \ref{tm:luzinplus}, 
\ref{tm:betaomega}, \ref{tm:p2}, \ref{tm:shz}, \ref{cor:irred} ) 
we don't need the full power of $\ccc {\kappa};$.
\begin{definition}\label{pr:c}
For $\kappa=\cf(\kappa)>\omega$   
 principle $\ccc {\kappa};$  ($\cccmin {\kappa};$) 
is the following statement:
\newline
For every $T\subs {\omega}^{<{\omega}}$ and
$\acal\in\matr {\kappa};$ we have
 (1) or (2) below:
\begin{enumerate}\arablabel
\item  \label{c1} there is a stationary (cofinal)
set $S\subs {\kappa}$ such that
if  $t\in T$ and  
$s\in \injs S;|t|;$ then 
$$
\aint s;t;\ne\empt,
$$
\item  \label{c2} there are  $t\in T$  and  stationary (cofinal) subsets
$D_0$, $D_1$, $\dots$, $D_{|t|-1}$ of ${\kappa}$ such that for every 
$s\in \injc D_0,\dots, D_{|t|-1};$ we have
$$
\aint s;t;=\empt.
$$
\end{enumerate}
\end{definition}

Next we formulate a dual version of principles 
$\ccc {\kappa};$ and $\cccmin {\kappa};$. Although  
we don't yet know any application of  principles $\ccci {\kappa};$
and $\cccimin {\kappa};$,  for the sake of
completeness  we include their definitions  here. 
Let us remark that we don't know whether 
$\ccc {\kappa};$ ($\cccmin {\kappa};$) implies 
$\ccci {\kappa};$  ($\cccimin {\kappa};$) or vice versa.
\begin{definition}\label{pr:ci}
If $\kappa=\cf(\kappa)>\omega$, then  
 principle $\ccci {\kappa};$  ($\cccimin {\kappa};$) 
is the following statement:
\newline For every $T\subs {\omega}^{<{\omega}}$ and
$\acal\in \matr {\kappa};$ we have
(1) or (2) below:
\begin{enumerate}\arablabel
\item  \label{cc1} there is a stationary (cofinal) 
set $S\subs {\kappa}$ such that
for each $t\in T$ and  
$s \in\injs S;|t|;$
$$
|\aint s;t;|<{\omega},
$$
\item \label{cc2} there are  $t\in T$  and  stationary (cofinal) subsets
$D_0$, $D_1$, $\dots$, $D_{|t|-1}$ of ${\kappa}$ such that for every
$s\in \injc  D_0,\dots, D_{|t|-1};$ we have
$$
|\aint s;t;|={\omega}.
$$
\end{enumerate}
\end{definition}

Let us remark that in the ``plain'' dual of principle $\ccc {\kappa};$ 
we should have 
 $|\aint s;t;|=\empt$ in \ref{pr:ci}\ref{cc1} and
$|\aint s;t;|\ne\empt$ in \ref{pr:ci}\ref{cc2}, but this 
``principle'' is   easily provable in ZFC.

The principles $\diamin {\kappa};$ and $\dia {\kappa};$ 
that we introduce next  easily follow from 
$\cccmin {\kappa};$ and $\ccc {\kappa};$, respectively, but as
their formulation is much simpler, we thought it to be worth while
to have them as separate principles. We first give two auxiliary
definitions.

\begin{definition}\label{pr:d}
If $\acal=\<\aaa i;{\alpha};:{\alpha}<{\kappa},i<{\omega}\>\in\matr {\kappa};$,
then we set  
$$
\lsup \acal;{\omega};=
\{Y\subs {\omega}:|\{{\alpha}<{\kappa}:\exists i<{\omega}\ 
\aaa i;{\alpha};\subs Y\}|= {\kappa}\}
$$
and
$$
\lsups \acal;{\omega};=
\{Y\subs {\omega}:\{{\alpha}<{\kappa}:\exists i<{\omega}\ 
\aaa i;{\alpha};\subs Y\}\text{ is stationary in }{\kappa}\}.
$$
\end{definition}

\begin{definition}\label{df:omega-adic}
A matrix $\acal=\<\aaa i;{\alpha};:{\alpha}\in S,i\in{\omega}\>$ is called
{\em ${\omega}$-adic} if for each 
 $t\in {\omega}^{<{\omega}}$ and  $s\in \injs S;|t|;$ 
we have $\aint s;t;\ne\empt$.
\end{definition}

Now we can formulate $\dia {\kappa};$ ($\diamin {\kappa};$)
as follows.

\begin{definition}\label{df:p-d}
For ${\kappa}=\cf({\kappa})>{\omega}$  principle 
 $\dia {\kappa};$ ($\diamin {\kappa};$ ) is the following statement: 
\newline If $\acal\in\matr {\kappa};$ and $\lsups \acal;{\omega};$
 ($\lsup \acal;{\omega};$) is centered then there is 
 a stationary (cofinal) set $S\subs {\kappa}$ such that 
$\acal\rest S$ is ${\omega}$-adic. 
\end{definition}

\begin{theorem}\label{tm:c-d}
$\ccc {\kappa};$ ($\cccmin {\kappa};$) implies $\dia {\kappa};$
($\diamin {\kappa};$).
\end{theorem}

\proof
We give  the proof only for $\dia {\kappa};$ because
the same argument works for $\diamin {\kappa};$. 

Let $\acal\in\matr {\kappa};$ and put $T={\omega}^{<{\omega}}$. 
By $\ccc {\kappa};$ either \ref{pr:c}(1) 
or \ref{pr:c}(2) holds. 

If $S\subs {\kappa}$ witnesses \ref{pr:c}(1)
for our $T$ then 
$\acal\rest S
$
is clearly ${\omega}$-adic. So it is enough to show that \ref{pr:c}(2) can not
hold. 

Assume, on the contrary, that there are  $t\in T={\omega}^{<{\omega}}$ and
 stationary subsets $D_0$, $D_1$, $\dots$, $D_{|t|-1}$ of ${\kappa}$
such that for each $s\in \injc D_0,\dots,D_{|t|-1};$ we have
\begin{equation}
\tag{+}\aint s;t;=\empt.
\end{equation}
We can obviously assume that the sets $D_i$ are pairwise disjoint.
Let $X_i=\bigcup\{\aaa t(i);{\delta};:{\delta}\in D_i\}$ 
for every $i<|t|$.
Then clearly $X_i\in\lsups \acal;{\omega};$ for every 
$i<|t|$, 
while (+) implies 
$\bigcap\limits_{i<k}X_i=\empt$, contradicting that 
$\lsups \acal;{\omega};$ is centered.  
\eproof

\begin{definition}\label{pr:e}
If ${\kappa}=\cf({\kappa})>{\omega}$, then 
 principle $\eee {\kappa};$ ($\eeemin {\kappa};$) is the following statement:
\newline
For every  $T\subs {\omega}^{<{\omega}}$ and
$\acal\in\matr {\kappa};$
(1) or (2) below holds:
\begin{enumerate}\arablabel
\item \label{f1} there is a stationary (cofinal) 
set $S\subs {\kappa}$ such that
$$
|\{\aint s;t;:t\in T\text{ and }s\in\injs S;|t|;\}|\le\omega.
$$
\item  \label{f2} there are  $t\in T$  and  stationary (cofinal) subsets
$D_0$, $D_1$, $\dots$, $D_{|t|-1}$ of ${\kappa}$ such that 
if $s_0,s_1\in \injc D_0,\dots, D_{|t|-1};$ with 
$s_0(i)\ne s_1(i)$ for each $i<|t|$ then  we have
$$
\aint s_0;t;\ne\aint s_1;t;.
$$
Clearly, if $t$, $D_0$, $\dots$,$D_{|t|-1}$ satisfy (2) then 
$$
|\{\aint s;t;:s\in\injc D_0,\dots,D_{|t|-1};\}|=\kappa.
$$
\end{enumerate}
\end{definition}

There is a surprising connection between these principles and the dual
versions $\ccci {\kappa};$ ($\cccimin {\kappa};$) of $\ccc {\kappa};$
($\cccmin {\kappa};$), respectively.


\begin{theorem}
$\eee {\kappa};$ $(\eeemin {\kappa};)$ implies 
$\ccci {\kappa};$ $(\cccimin {\kappa};)$.
\end{theorem}

\proof
Let $\acal\in\matr {\kappa};$ and $T\subs {\omega}^{<{\omega}}$ and
apply $\eee {\kappa};$ to $\acal$ and $T$. Assume first  
that there is a stationary set $S\subs {\kappa}$ such that the family
$$
\ical=\{\aint s;t;:t\in T\text{ and }s\in\injs S;|t|;\}
$$
is countable. 

Now  for  $t\in T$, $i<|t|$ and 
$I\in\ical\cap \br {\omega};{\omega};$
set   
$$
D(I,t,i)=\{{\alpha}\in S:\aaa t(i);{\alpha};\sups I\}.$$ 
If for some $t\in T$ and $I\in \ical\cap \br {\omega};{\omega};$
the set $D(I,t,i)$ is stationary for each $i<|t|$ then this $t$ and the sets 
$D(I,t,0)$, $\dots$, $D(I,t,|t|-1)$ witness  \ref{pr:ci}\ref{cc2}.

So we can assume that for all $t\in T$ and 
$I\in \ical\cap \br {\omega};{\omega};$ the set
$$
b(I,t)=\{i<|t|:D(I,t,i)\text{ is non-stationary in }{\kappa}\}
$$ 
is not empty. Then the set 
$$
D=\bigcup\{D(I,t,i):I\in\ical\cap \br {\omega};{\omega};, t\in T,
i\in b(I,t)\}
$$
is not stationary and so $S'=S\setm D$ is stationary. 
We claim that $S'$ witnesses \ref{pr:ci}\ref{cc1}.
 Assume on the contrary that $t\in T$,  $s\in \injs S';|t|;$ and 
 $I=\aint s;t;$ is infinite. Then $I\in\ical\cap \br {\omega};{\omega};$ and 
$s(i)\in D(I,t,i)$ for each $i<|t|$. Since $s(i)\notin D$ it follows
that  $D(I,t,i)$ is stationary for each $i<|t|$, that is,
$b(I,t)=\empt$, which is a contradiction. 

Assume now that there are  $t\in T$  and  stationary subsets
$D_0$, $D_1$, $\dots$, $D_{|t|-1}$ of ${\kappa}$ such that 
if $s_0,s_1\in \injc D_0,\dots, D_{|t|-1};$ with 
$s_0(i)\ne s_1(i)$ for each $i<|t|$ then  
$$
\aint s_0;t;\ne\aint s_1;t;.
$$
We show that in this case again \ref{pr:ci}\ref{cc2} holds.
Indeed,  for each $I\in \br {\omega};<{\omega};$ pick  
$s_I\in \injc D_0,\dots, D_{|t|-1};$ such that 
$\aint s_I;t;=I$ provided that there is such an $s$.
Let $R=\bigcup\{s_I(i):I\in \br {\omega};<{\omega};, i<|t|\}$
and $D'_i=D_i\setm R$ for $i<|t|$.
Now if $s\in \injc D'_0,\dots, D'_{|t|-1};$ then for any 
$I\in \br {\omega};<{\omega};$ we have $s_I(i)\ne s(i)$ for
each $i<|t|$, hence  
$I=\aint s_I;t;\ne \aint s;t;$.
As $I$ was an arbitrary element of $\br {\omega};<{\omega};$ 
we  conclude that  $|\aint s;t;|={\omega}.$
\eproof

If ${\kappa}=\cf({\kappa})>c$ then 
$\ccc {\kappa};$ and
$\eee {\kappa};$ are trivially valid.
Indeed given  $\acal\in\matr {\kappa};$ and $T\subs  {\omega}^{<{\omega}}$
there is a stationary set $S\subs {\kappa}$ such that for any
${\alpha},{\beta}\in S$ and $n\in {\omega}$ we have
$\aaa n;{\alpha};=\aaa n;{\beta};$. Then $S$ witnesses \ref{pr:e}(1)
and so principle $\eee {\kappa};$ holds. 
If $S$ does not witness \ref{pr:c}(1) 
then for some $t=\nok\in T$
we have $\bigcap\limits_{i<k}A_{n_i}=\empt$. Thus  $t$ and
$D_0=D_1=\dots=D_{k-1}=S$ witness \ref{pr:c}(2).

As was mentioned in section \ref{sc:int}, our principles are of
interest only for $ {\kappa} > \oo $. In fact, for ${\kappa}=\oo$,
they are all false! 


To see that $\diamin \oo;$ $($so also $\cccmin \oo;$$)$ is false  we
may recall that
in \cite{JNSS} we have constructed, in ZFC, a separable, first countable 
$P_{<{\omega}}$ space $X$ of size $\oo$. 
(A Hausdorff space $X$ is called $P_{<{\omega}}$ 
if the intersection of finitely many uncountable open subsets of $X$ is
always non-empty.) We can assume that the  underlying set of $X$ is $\oo$
and ${\omega}$ is dense in $X$. For each ${\alpha}<\oo$ let 
$\{U({\alpha},n):n<{\omega}\}$ be a neighbourhood base of ${\alpha}$ in $X$. 
Now consider the $\oo\times {\omega}$-matrix
$$
\acal=\<U({\alpha},n)\cap {\omega}:{\alpha}<\oo, n<{\omega}\>.
$$
Then $B\in\lsup \acal;{\omega};$ if and only if there is an uncountable
open set $U\subs X$ such that $U\cap {\omega}\subs B$.
Since $X$ is a $P_{<{\omega}}$ space it follows that 
$\lsup \acal;{\omega};$ is centered. 
But the space $X$ is Hausdorff, so there is not even a two
element subset $S$ of $X$ such that 
$\acal\rest S$ is ${\omega}$-adic.  


To show that $\eeemin \oo;$ $($and so also $\cccimin \oo;$$)$ is false 
we need the following observation.

\begin{theorem}\label{tm:centered-notlinked}
There is a subfamily
$\acal=\{A_{\alpha}:{\alpha}<\oo\}$ of  $\br {\omega};{\omega};$ 
such that for any $n\in {\omega}$ and $I_0,\dots, I_{n-1}\in \br \oo;\oo;$
there are ${\gamma}_i, {\delta}_i\in I_i$ for $i<n$
with $\bigcap\{A_{{\gamma}_i}:i<n\}$ is infinite but 
$A_{{\delta}_i}\cap A_{{\delta}_j}$ is finite for any $i<j<n$.
\end{theorem}

\proof
The proof is based on two lemmas which are   
 probably well-known.
\begin{lemma}\label{lm:f-sier}
There is a function $f:\br\oo;2;\to 2$ such that 
for any $n\in {\omega}$ and $I_0,\dots, I_{n-1}\in \br \oo;\oo;$
there are ${\gamma}_i, {\delta}_i\in I_i$ for $i<n$
such that $f({\gamma}_i,{\gamma}_j)=0$ and 
$f({\delta}_i,{\delta}_j)=1$  for each $i\ne j<n$.
\end{lemma}

\proof
We show that the Sierpienski coloring has this property. 
So let $\{r_{\alpha}:{\alpha}<\oo\}$ be pairwise different real numbers
and for ${\alpha}<{\beta}<\oo$ put $f({\alpha},{\beta})=0$ iff 
$r_{\alpha}<r_{\beta}$.
Given  $n\in {\omega}$ and $I_0,\dots, I_{n-1}\in \br \oo;\oo;$
 let  $x_i$ be a complete  accumulation point of  
$A_i=\{r_{\alpha}:{\alpha}\in I_i\}$ with $x_i\ne x_j$ 
for $i<j<n$. We may assume that 
$x_0<x_1<\dots<x_{n-1}$. So there are
$I'_i\in\br I_i;\oo;$ for $i<n$ such that if 
$A'_i=\{r_{\alpha}:{\alpha}\in I'_i\}$ then $A'_i<_R A'_j$
whenever $i<j<n$. Now pick first ${\gamma}_{0}\in I'_0$, then
${\gamma}_{1}\in I'_1\setm ({\gamma}_{0}+1)$, then
${\gamma}_{2}\in I'_2\setm ({\gamma}_{1}+1)$ and so on.
Then we have $f({\gamma}_{i},{\gamma}_{j})=0$ for each $i<j<n$.
Next we pick ${\delta}_{n-1}\in I'_{n-1}$, then
${\delta}_{n-2}\in I'_{n-2}\setm {\delta}_{n-1}$, then
${\delta}_{n-3}\in I'_{n-3}\setm {\delta}_{n-2}$ and so on.
Then we have $f({\delta}_{i},{\delta}_{j})=1$ for  $i<j<n$.
\eproof

\def\aij#1;#2;{A[#1,#2]}
\def\aii#1;{A[#1]}

A family $\ical\subs \pcal(\oo)$ is called {\em downwards closed}
if $\pcal(I)\subs \ical$  for each $I\in\ical$.
Given a family $\acal=\{A_{\alpha}:{\alpha}<\oo\}\subs \pcal({\omega})$  
and $I,J\in \br \oo;<{\omega};$,  let 
$$
\aij I;J;=\bigcap\{A_{\alpha}: {\alpha}\in I\} \setm
\bigcup\{A_{\beta}:  {\beta}\in J\}.$$
Put $\aii I;=\aij I;\empt;$.

Clearly, if $\acal\subs \pcal({\omega})$ then
$\{I\in \br {\omega};<{\omega};: |\aii I;|={\omega}\}$ is downward
closed. Our next result is a converse of this statement.

\begin{lemma}\label{lm:repr}
If $\ical\subs \br \oo;<{\omega};$ is downwards closed then there
is a family 
$\acal=\{A_{\alpha}:{\alpha}<\oo\}$ of subsets of ${\omega}$ 
such that 
\begin{equation}
\tag{\dag}
\ical=\{I\in \br\oo;<{\omega};:\aii I;\text{ is infinite}\}.
\end{equation}
\end{lemma}

\proof
For ${\alpha}<\oo$
write $\ical_{\alpha}=\ical\cap \br {\alpha}+1;<{\omega};$.

 We will define 
 $A_{\alpha}\subs {\omega}$ by induction on ${\alpha}<\oo$ so as to
satisfy the following inductive hypotheses $(\ddag)_{\alpha}$ which is
stronger than $(\dag)$ restricted to ${\alpha}$:
\begin{equation}
\tag{$\ddag_{\alpha}$} 
\ical_{\alpha}=\{I\in \br {\alpha}+1;<{\omega};:
\forall J\in \br {\alpha}\setm I;<{\omega};\ \aij I;J;\text{ is infinite.}\}
\end{equation}
 
The induction uses the following elementary fact. 

\begin{fact}\label{fact}
If $\bcal$  and $\dcal$ are countable subfamilies 
 of $\br {\omega};{\omega};$ such that no element of $\bcal$
is covered by the union of finitely many elements of $\dcal$ 
then there is a set $X\subs {\omega}$ such that 
\begin{enumerate}\rlabel
\item $B\cap X$ is infinite for each $B\in \bcal$,
\item $B\setm X$ is infinite for each $B\in \bcal$,
\item $D\cap X$ is finite for each $D\in \dcal$.
\end{enumerate}
\end{fact}

Now, if $A_{\beta}$ has been defined and $(\ddag_{\beta})$ holds for 
all ${\beta}<{\alpha}$ then let 
$\bcal=\{\aij I;J;:I\in \ical_{\alpha}, J\in \br {\alpha};<{\omega};
\land I\cap J=\empt\}$,
$\dcal=\{\aij I;J;:I\in \br {\alpha};<{\omega};\setm \ical_{\alpha}, 
J\in \br {\alpha};<{\omega};\land
I\cap J=\empt\}$ and apply  fact \ref{fact} to get $A_{\alpha}$.
It is easy to check that $(\ddag_{\alpha})$ will be satisfied.
 \eproof

To get a family $\acal$  satisfying the requirements of theorem 
\ref{tm:centered-notlinked} take the function $f$ given by lemma 
\ref{lm:f-sier}  and apply  lemma \ref{lm:repr} to  
$\ical=\{I\in \br \oo;<{\omega};:f''I\subs\{0\}\}$.
\eproof

Theorem \ref{tm:centered-notlinked} yields immediately 
the following corollary:

\begin{corollary}\label{cor:centered-notlinked}
There is a family 
$\acal=\{A_{\alpha}:{\alpha}<{\omega}_1\}$
such that 
$$
\lsup \acal;{\omega};=
\{Y\subs {\omega}:|\{{\alpha}<{\omega}_1:A_{\alpha}\subs Y\}|={\omega}_1\}
$$
is centered but for  no  
$S\in \br \oo;\oo;$ is $\{A_{\alpha}:{\alpha}\in S\}$ 
 linked.
\end{corollary}

Let us remark that if $\acal$ is almost disjoint then 
 $\lsup \acal;{\omega};$ is centered if and only if
$\acal$ is a strong Luzin gap (i.e. there is no partition of
${\omega}$ into finitely many pieces such that
each piece is almost disjoint to uncountably many elements of $\acal$).
If $\text{MA}_{\oo}$ holds, then there is no strong Luzin gap
(see \cite[Theorem 3.2]{JNSS}), so  in ZFC one can not construct an 
almost disjoint family $\acal$ satisfying the requirements 
of corollary \ref{cor:centered-notlinked}.

The family $\acal$ of \ref{cor:centered-notlinked} can be used to give
counterexamples to both $\diamin \oo;$ and $\cccmin \oo;$, in fact via
th same matrix in $\matr \oo;$.

\begin{corollary}
$\cccimin \oo;$ (and so $\eeemin \oo;$ too) and $\diamin \oo;$ 
are both false.
\end{corollary}

\proof
Consider the family $\acal=\{A_{\alpha}:{\alpha}<\oo\}$ given
by  \ref{cor:centered-notlinked}.   
Put  $T={\omega}^{<{\omega}}$ and
$\aaa i;{\alpha};=A_{\alpha}$  for each ${\alpha}<\oo$ and $i<{\omega}$.
Then neither \ref{pr:ci}\ref{cc1} nor \ref{pr:ci}\ref{cc2} can hold for 
$\acal=\<\aaa i;{\alpha};:{\alpha}<\oo,i<{\omega}\>$ and $T$.
Moreover, the matrix $\acal$ clearly contradicts $\diamin \oo;$. 
\eproof

\section{Consistency of the principles in the Cohen model.}

A cardinal ${\kappa}$ is  
{\em ${\omega}$-inaccessible} if
 ${\lambda}^{\omega}<{\kappa}$ holds for each
${\lambda}<{\kappa}$.
Given any set $I$ we denote by $\coh I;$ the poset $\fn(I,2,{\omega})$,
i.e. the standard one adding ${\kappa}$ Cohen reals.

In this section we prove that if 
${\kappa}$ is a regular ${\omega}$-inaccessible cardinal in some
ground model $V$ and we add any number of  Cohen reals to $V$
 then in the  extension
the principles $\ccc {\kappa};$, $\ccci {\kappa};$ and
$\eee {\kappa};$  are all satisfied. 
As we remarked in section \ref{sc:principles}   above  
the  case  ${\kappa}>{\lambda}$  is trivial, while  the 
 case ${\kappa}<{\lambda}$  can be reduced to the  
case ${\kappa}={\lambda}$.

Since the proof of the latter is  long and technical, 
we first sketch the main
idea. So let us be given  a matrix $\acal\in\matr {\kappa};$
and a set $T\subs {\omega}^{<{\omega}}$ in $V[G]$, where
$G$ is $\coh {\kappa};$-generic over $V$.
In the first part of the proof we find a set  $I\in\br {\kappa};{\omega};$
and a stationary
set $S\subs {\kappa}$ such that in $V[G\rest I]$ the sequences 
$\<\aaa i;{\alpha};:i<{\omega}\>$ for ${\alpha}\in S$ have 
also pairwise isomorphic names
with disjoint supports (contained in ${\kappa}\setm I$).
This reduction,  carried out in   lemma \ref{lm:reduction},
 will be the place where we  use that ${\kappa}$
is regular and ${\omega}$-inaccessible in V.
In the second part of the proof,  using slightly different
arguments for $\ccc {\kappa};$ and for $\eee {\kappa};$,  we show 
that if some $\acal\in\matr S;$ has names with these properties then
either  $S$ witnesses \ref{pr:c}(1) 
(or \ref{pr:e}(1), respectively)
  or some stationary sets   $D_i\subs S$  witness 
\ref{pr:c}(2) (or \ref{pr:e}(2), respectively).
In this second step we don't use that
${\kappa}$ is ${\omega}$-inaccessible or  regular.

In our forcing  arguments  
we  follow the notation of Kunen \cite{Ku}.
Let us first recall  definition \cite[5.11]{Ku}. 

\begin{definition}
A $\coh I;$-name $\bdot$ of a subset of some ordinal ${\mu}$ is called 
{\em nice} if for each ${\nu}<{\mu}$ there is an antichain 
$B_{\nu}\subs \coh I;$ such that 
$$
\bdot=\{\<p,\hat{{\nu}}\>:{\nu}\in {\mu}\land p\in B_{\nu}\}=
\bigcup\{B_{\nu}\times\{\hat{{\nu}}\}:{\nu}\in {\mu}\}.
$$
We let $\supp(\bdot)=\bigcup\{\dom(p):p\in\bigcup\limits_{{\nu}<{\mu}}B_{\nu}\}$.
\end{definition}

It is well-known (see e.g. lemma \cite[5.12]{Ku}) 
that every set of ordinals in  $V[G]$ has a nice name in $V$.

If ${\varphi}$ is a bijection between two sets $I$ and $J$ then
${\varphi}$ lifts to a natural isomorphism between $\coh I;$ and
$\coh J;$, which will be also denoted by ${\varphi}$, as follows: 
for $p\in \coh I;$ let $\dom({\varphi}(p))={\varphi}''\dom(p)$ and
${\varphi}(p)({\varphi}({\xi}))=p({\xi})$. Moreover ${\varphi}$ also
generates a bijection between the nice $\coh I;$-names and the nice
$\coh J;$-names (see \cite[7.12]{Ku}): 
if $\bdot$ is a nice $\coh I;$-name then let
${\varphi}(\bdot)=\{\<{\varphi}(p),\hat{{\xi}}\>:\<p,\hat{{\xi}}\>\in\bdot\}$.
\def\fij#1;#2;{{\varphi}_{#1,#2}}
If $I$ and $J$ are sets of ordinals with the same order type then
  $\fij I;J;$ is the natural order-preserving
bijection from $I$ onto $J$.

\begin{definition}\label{df:twin-names}
Assume $I,J\subs {\kappa}$,  moreover $\adot_i$ and $\bdot_i$ are nice 
$\coh {\kappa};$-names of subsets of ${\omega}$ for $i<{\omega}$,
such that  $\supp(\adot_i)\subs I$ and $\supp(\bdot_i)\subs J$. 
We say that the structures of names $\<I,\adot_i:i<{\omega}\>$ and  $\<J,\bdot_i:i<{\omega}\>$ are {\em twins} 
if $I$ and $J$ have the same order type and
\begin{enumerate}\arablabel
\item $\fij I;J;$ is the identity on $I\cap J$,
\item  $\fij I;J;(\adot_i)=\bdot_i$ for each $i<{\omega}$.
\end{enumerate}
\end{definition}


\begin{definition}\label{df:projection}
Assume that $I\subs {\kappa}$, $G$ is a $\coh {\kappa};$-generic filter 
over $V$  and $H=G\rest I$.
If $\bdot$ is a   nice $\coh {\kappa};$-name  of a subset of some ordinal 
${\mu}$ we define
 in $V[H]$  the $\coh {\kappa}\setm I;$ name $\pii H;{\kappa};I; (\bdot)$ as follows:
$$
\pii H;{\kappa};I;(\bdot)=\{\<p\rest {\kappa}\setm I,\hat{{\nu}}\>:
\<p,\hat{{\nu}}\>\in\bdot\land
p\rest I\in H\}.
$$  
\end{definition}

\begin{lemma}\label{lm:projection}
$\pii H;{\kappa};I;(\bdot)$ is a nice $\coh {\kappa}\setm I;$-name in $V[H]$
and $\supp(\pii H;{\kappa};I;(\bdot))\subs\supp(\bdot)\setm I$, moreover
$$\ig G\rest ({\kappa}\setm I);{\pii H;{\kappa};I;(\bdot)};=\ig G;{\bdot};.
$$
\end{lemma}

\proof
Straightforward from the construction.
\eproof

\begin{definition}\label{df:nice-matrix}
Assume that $S\subs {\kappa}$.
A matrix $\bcaldot=\<\bbbdot i;{\alpha};:{\alpha}\in S,i<{\omega}\>$
of nice $\coh {\kappa};$-names of subsets of ${\omega}$
is called a {\em nice $S$-matrix} if 
 conditions (i) and (ii) below hold: 
\begin{enumerate}\rlabel
\item putting $J_{\alpha}=\bigcup_{i<{\omega}}\supp(\bbbdot i;{\alpha};)$ 
the sets $\{J_{\alpha}:{\alpha}\in S\}$ are pairwise disjoint,
\item  the structures of names 
$\{\<J_{\alpha},\bbbdot i;{\alpha};:i<{\omega}\>:{\alpha}\in S\}$
 are pairwise twins. 
\end{enumerate}
We denote by $\matrn S;$ the family of nice $S$-matrices.
\end{definition}

\begin{lemma} {\em \bf(Reduction lemma)}
\label{lm:reduction}
Assume that ${\kappa}$ is a regular, ${\omega}$-inaccessible cardinal, 
$G$ is $\coh {\kappa};$-generic over $V$ and 
$\acal\in\matr {\kappa};$ in $V[G]$.
Then there are a countable set $I\subs {\kappa}$ and  a 
 stationary set $S\subs {\kappa}$ in $V$ such that,   in
$V[G\rest I]$, there is  $\bcal\in \matrn S;$ satisfying    
$V[G]\models$ ``$\aaa i;{\alpha};=
\ig G\rest ({\kappa}\setm I);{\bbbdot i;{\alpha};};$''
 for each ${\alpha}\in S$ and $i\in {\omega}$.
\end{lemma}

\proof
Assume that 
$$
1_{\coh {\kappa};}\force\text{
``$\dot\acal=\<\aaadot i;{\alpha};:{\alpha}<{\kappa},i<{\omega}\>
\in\matr {\kappa};$.'' 
}
$$
We can assume that all the names $\aaadot {\alpha};i;$ are nice.
Let $I_{\alpha}=\bigcup_{i<{\omega}}\supp(\aaadot i;{\alpha};)$. 

We need a strong version  of Erd\H os-Rado $\Delta$-system theorem 
 saying  that there is a stationary set
$T\subs {\kappa}$ such that $\{I_{\alpha}:{\alpha}\in T\}$
forms a $\Delta$-system with 
some kernel $I$, moreover $\sup I<\min I_{\alpha}\setm I$ for each 
${\alpha}\in T$.
Although this statement is well-known we present a proof because
we could not find any reference to it.

\begin{ertheorem}\label{lm:erdosrado}
If ${\kappa}$ is an ${\omega}$-inaccessible regular cardinal and $\xcal=\{X_{\alpha}:{\alpha}<{\kappa}\}$ is a family 
of  countable sets then there is a stationary set $I\subs {\kappa}$ such that 
$\{X_{\alpha}:{\alpha}\in I\}$ forms a $\Delta$-system.
\end{ertheorem}

\begin{pf}
Since $|\cup\acal|\le{\kappa}$ we can assume that 
$X_{\alpha}\subs {\kappa}$. Let $J=\{{\alpha}<{\kappa}:\cf({\alpha})=\oo\}$.
Define the function $f:J\to {\kappa}$ by the stipulation
$f({\alpha})=\sup (X_{\alpha}\cap A)$. Since  $X_{\alpha}$ is countable 
and $\cf({\alpha})=\oo$ we have $f({\alpha})<{\alpha}$, i.e. the function
$f$ is regressive on the stationary set $J$. So by the Fodor lemma,
$f$ is constant on a stationary set $K\subs J$. 
Say $f'' K=\{{\nu}\}$.
For ${\alpha}\in K$ let
$h({\alpha})=X_{\alpha}\cap {\nu}$. Since ${\nu}<{\kappa}$ it follows that 
 the range of $h$ is of size $|{\nu}|^{\omega}<{\kappa}$. 
But $K$ is stationary, so there is a stationary $M\subs K$ such that $h$ is constant on $M$, say $h''M=\{A\}$.

For ${\alpha}\in {\kappa}$ let 
$g({\alpha})=\sup X_{\alpha}$. Then the set  
$$
C=\{{\beta}<{\kappa}: g({\gamma})<{\beta}\text{ for each }{\gamma}<{\beta} \}
$$
is club in ${\kappa}$. Let $I=M\cap C$.  
We show that $\{X_{\alpha}:{\alpha}\in I\}$ forms a $\Delta$-system
with kernel $A$. Let ${\alpha},{\beta}\in I$, ${\alpha}<{\beta}$.
Then $X_{\alpha}\cap X_{\beta}=(X_{\alpha}\cap (X_{\beta}\cap {\beta}))\cup
(X_{\alpha}\cap(X_{\beta}\setm {\beta})$.

But $(X_{\beta}\cap {\beta})=A$ and so 
$(X_{\alpha}\cap (X_{\beta}\cap {\beta}))=A$. Since ${\beta}\in C$
it follows that $g({\alpha})<{\beta}$, i.e. $X_{\alpha}\subs {\beta}$
and so $(X_{\alpha}\cap(X_{\beta}\setm {\beta})=\empt$.
Putting together these two equations we obtain
$X_{\alpha}\cap X_{\beta}=A$ which was to be proved.
\end{pf}

Since $2^{\omega}<{\kappa}=\cf({\kappa})$ and there are only 
$2^{\omega}$ different isomorphism types of structures 
of names there is a stationary set $S\subs T$ such that
the structures of names 
$\{\<I_{\alpha},\aaadot i;{\alpha};:i<{\omega}\>:{\alpha}\in S\}$ 
are pairwise twins.

From now on we work in $V[G\rest I]$. Let 
$\bbbdot i;{\alpha};=\pii G\rest I;{\kappa};I;(\aaadot i;{\alpha};)$ for ${\alpha}\in S$ and $i\in {\omega}$. Then 
  $\supp(\bbbdot i;{\alpha};)\subs J_{\alpha}=I_{\alpha}\setm I$ 
and the structures
of names $\<J_{\alpha},\bbbdot i;{\alpha};:i<{\omega}\>$ are pairwise
twins by lemma \ref{lm:projection} above.  Thus  
$\bcal=\<\bbbdot i;{\alpha};:{\alpha}\in S,i<{\omega}\>\in\matrn S;$.
\eproof

\def\dotb#1;{\bdot_{#1}}
\def\seqn#1;{\scal(#1)}
\def\phis#1;#2;{{\varphi}({\dotb {#1(0)};},\dots,{\dotb {#1(#2)};},\hat{Z})}
\def\phisev#1;#2;{{\varphi}({B_{#1(0)}},\dots,{B_{#1(#2)}},Z)}

\begin{definition}\label{df:nice-sequence}
Assume that $S\subs {\kappa}$.
A sequence  $\bcaldot=\<\<J_{\alpha},\dotb {\alpha};\>:{\alpha}\in S\>$
is called a {\em nice $S$-sequence} if 
 conditions (i) and (ii) below hold: 
\begin{enumerate}\rlabel
\item $J_{\alpha}\in \br {\kappa};{\omega};$, $\dotb {\alpha};$
is a nice $\coh J_{\alpha};$-name,  and $J_{\alpha}$
for ${\alpha}\in S$ are pairwise disjoint,
\item  the structures of names 
$\<J_{\alpha},\dotb {\alpha};\>$
for ${\alpha}\in S$ are pairwise twins. 
\end{enumerate}
We denote by $\seqn S;$ the family of nice $S$-sequences.
\end{definition}

\begin{lemma}{\em \bf (Homogeneity lemma)}
\label{lm:homog}
Assume that $S\subs {\kappa}$ and   
$\bcaldot=\<\<J_{\alpha},\dotb {\alpha};\>:{\alpha}\in S\>$ is a 
 {\em nice $S$-sequence}. If ${\varphi}(x_0, x_1,\dots, x_{n-1},z)$ 
is a formula with free variables $x_0, x_1,\dots, x_{n-1},z$ and
$Z$ is an element of the ground model, then
$(1)$ or $(2)$
below holds:
\begin{enumerate}\arablabel
\item  $1_{\coh {\kappa};}\force$ ``
${\varphi}(\dotb s(0);\dotb s(1);,\dots,\dotb s(n-1);,\hat{Z})$
for all  $s\in\injs S;n;$'', 
\item  for some $r\in \coh {\kappa};$  we have \newline 
$r\force$ ``
there are subsets $\dot D_0$, $\dot D_1$, $\dots$, $\dot D_{n-1}$ of 
$S$ such that
\begin{enumerate}\alabel
\item    for each $i<n$ and  $A\in \br S;{\omega};\cap V$ 
we have $\dot D_i\cap A\ne \empt$,
\item $\neg{\varphi}(\dotb s(0);\dotb s(1);,\dots,
\dotb s(n-1);,\hat{Z})$ for all  
$s\in\injc \dot D_0, \dot D_1, \dots, \dot D_{k-1};$''.
\end{enumerate} 
\end{enumerate}
\end{lemma}

\proof
Assume that (1) fails, that is, 
there are $p\in \coh {\kappa};$ and  
$s\in\injs S;k;$ such that 
$$
p\force \text{``}\neg \phis s;n-1;\text{''.}
$$
Let $J=\bigcup\limits_{i<k}J_{s(i)}$ and $p'=p\rest J$ and
$r=p\setm p'$.
Since the sets $J_{\alpha}$  are pairwise disjoint we can
assume that $\dom (r)\cap J_{\alpha}=\empt$ for each ${\alpha}\in S$.

For $\<{\alpha},{\beta}\>\in S^2$ we denote  by
$\fab$ the natural order preserving bijection between $J_{\alpha}$
and $J_{\beta}$. 
For ${\beta}\in S$ and $i<k$ let
 $p({\beta},i)={\varphi}_{s(i),{\beta}}(p\rest J_{s(i)})$.
For $i<k$ define the $\coh {\kappa};$-name $\dotd_i$ of a subset of $S$
as follows:
 $\dotd_i=\{\<p({\beta},i),\hat{{\beta}}\>:{\beta}\in S\}$.
Then 
$$
V[G]\models \text{``}D_i=\ig G;\dotd_i;=
\{{\beta}\in S:p({\beta},i)\in G\}\text{''},
$$
where $G$ is  $\coh {\kappa};$-generic over $V$. Since
the supports of $p({\beta},i)$ for ${\beta}\in S$  are pairwise disjoint 
a standard density argument gives that $D_i\cap A\ne\empt$ 
whenever  $A\in \br S;{\omega};\cap V$, hence (a) holds.  

To show (b) assume that $r\in G$ and 
$$
V[G]\models u\in\injc  D_0,\dots,D_{k-1};.
$$
Since $u$ is finite we have $u\in V$.
Let $J^*=\bigcup\limits_{i<k}J_{u(i)}$ and
${\psi}=\bigcup\limits_{i<k}{\varphi}_{u(i),s(i)}$.
Then ${\psi}$ is a bijection between $J^*$ and $J$ and so it extends
to  isomorphisms between $\coh J^*;$ and $\coh J;$, and between the families
of nice $\coh J^*;$-names and of nice $\coh J;$-names.
Let $\Psi$ be the natural extension of $\psi$ to a permutation of ${\kappa}$:
$$\Psi({\nu})= \left\{ 
\begin{array}{ll}
\psi({\nu})&
\mbox{if ${\nu}\in J^*$},\\
\psi^{-1}({\nu})&
\mbox{if ${\nu}\in J$,}\\
{\nu}&\mbox{if ${\nu}\in{\kappa}\setm (J\cup J^*).$}
\end{array}
\right.
$$
Then $\Psi$  extends
to an automorphism of  $\coh {\kappa};$, and also to an  automorphism of 
 nice $\coh {\kappa};$-names. Clearly if $q\in \coh J^*;$ 
and $\bdot$ is a nice  $\coh J^*;$-name then 
$\psi(q)=\Psi(q)$ and $\psi(\bdot)=\Psi(\bdot)$.
Observe that $\Psi(r)=r$ and $\Psi(\hat{Z})=\hat{Z}$.

Let $G^*={\Psi}''G$. 
Then $G^*$ is also a $\coh {\kappa};$-generic filter over  $V$ and
since ${\Psi}(\dotb u(i);)=\dotb s(i);$ it follows that
\begin{equation}
\tag{$\bullet$}
\ig G;{\dotb u(i);};=\ig G^*;{\dotb s(i);};.
\end{equation}
But $p(u(i),i)\in G$, so 
$p\rest J_{s(i)}={\psi}(p(u(i),i))\in G^*$.
Thus $p=r\cup\bigcup\limits_{i<k}p\rest J_{s(i)}\in G^*$ as well.
Since $p\force$ $\neg\phis s;n-1;$ and so
$V[G^*]\models$ 
``$\neg\phisev s;n-1; $'', 
by $(\bullet)$ this implies 
$$
V[G]\models\text{``}\neg \phisev u;n-1;\text{''}
$$ 
which was to be proved.
\eproof

\begin{theorem}\label{tm:c}
If ${\kappa}$ is a regular, ${\omega}$-inaccessible cardinal 
then for each cardinal ${\lambda}$ we have 
$$
V^{\coh {\lambda};}\models\text{$\ccc {\kappa};$ and $\ccci {\kappa};$ hold.}
$$
\end{theorem}

\proof
We deal only with $\ccc {\kappa};$ because the same argument works for
$\ccci {\kappa};$. 
As we observed in section \ref{sc:principles} 
we can assume that $\kappa\le\lambda$.
First we investigate the case ${\lambda}={\kappa}$.

Assume that 
$$
1_{\coh {\kappa};}\force\text{
``$\dot\acal=\<\aaadot i;{\alpha};:{\alpha}<{\kappa};i<{\omega}\>
\in\matr {\kappa};$ and $T\subs {\omega}^{<{\omega}}$.'' 
}
$$
Applying the reduction lemma \ref{lm:reduction} and that $T$ is countable
we can find 
 a countable set $I\subs {\kappa}$ and 
 a stationary set $S\subs {\kappa}$ in $V$ and   
a nice $S$-matrix $\bcal$ in $V[G\rest I]$ 
such that   
$$
V[G]\models \text{``}\ig G;{\aaadot i;{\alpha};};=
\ig G\rest ({\kappa}\setm I);{\bbbdot i;{\alpha};};\text{''}
$$
 for  ${\alpha}\in S$ and $i\in {\omega}$, moreover
$T\in V[G\rest I]$.

We show that for each $q\in\coh {\kappa};$
there is a condition $r\le q$ in $\coh {\kappa};$ such that
$r\force $ ``\ref{pr:c}\ref{c1} or \ref{pr:c}\ref{c2} holds''.
Let $I'=I\cup \dom (q)$. 

For each $t\in T$ let ${\varphi}_t(x_0,\dots,x_{|t|-1})$ be the following
formula:
$$
{\varphi}(\<B_{0,k}:k<{\omega}\>,\dots,\<B_{|t|-1,k}:k<{\omega}\>)
\Longleftrightarrow \bigcap_{i<|t|} B_{i,t(i)}\ne\empt.
$$
Applying the homogeneity lemma \ref{lm:homog} to $V[G\rest I']$ 
as our ground model and to every ${\varphi}_t$ we get that either 
$q\force $ ``\ref{pr:c}\ref{c1} holds'' or 
$q\cup p\force$ ``\ref{pr:c}\ref{c2} holds''.
Let us remark that \ref{lm:homog}(2)(a) implies 
that  as $S$ is stationary, so is each $D_i$.

Thus  we have proved the theorem in the case ${\kappa}={\lambda}$.
If ${\lambda}>{\kappa}$ and
$\acal\in(\matr {\kappa};)^{V[G]}$, where
$G$ is $\coh {\lambda};$-generic over $V$, then 
there is $J\in \br {\lambda};{\kappa};$ 
such that $\acal\in  V[G\rest  J]$. 
The stationary sets that witness
\ref{pr:c}(1) or \ref{pr:c}(2) in $V[G\rest J]$ remain stationary in 
$V[G]$, and so we are done.
\eproof

\begin{theorem}\label{tm:e}
If ${\kappa}$ is a regular, ${\omega}$-inaccessible cardinal 
then for each cardinal ${\lambda}$ we have 
$$
V^{\coh {\lambda};}\models\text{$\eee {\kappa};$ holds.}
$$
\end{theorem}

\proof
As in  \ref{tm:c} the  important case  is when ${\lambda}={\kappa}$,
because the  case ${\lambda}<{\kappa}$  is trivial and the 
case ${\kappa}<{\lambda}$  can be reduced to the case ${\kappa}={\lambda}$.

So  assume that  
$$
1_{\coh {\kappa};}\force\text{
``$\dot\acal=\<\aaadot i;{\alpha};:{\alpha}<{\kappa},i<{\omega}\>
\in\matr {\kappa};$.'' 
}
$$

Applying lemma \ref{lm:reduction} we can find 
 a countable set $I\subs {\kappa}$, a 
 stationary set $S\subs {\kappa}$ in $V$ and  in
$V[G\rest I]$  a nice $S$-matrix $\bcal$ such that   
$V[G]\models$ ``$\aaa i;{\alpha};=
\ig G\rest ({\kappa}\setm I);{\bbbdot i;{\alpha};};$''
 for each ${\alpha}\in S$ and $i\in {\omega}$, moreover
$T\in V[G\rest I]$.

We need the following lemma that is probably well-known.

\begin{lemma}\label{lm:V}
If $H$ is  a $\coh {\kappa};$-generic filter  over $V$ and  
$I$, $J$ are disjoint
subsets of ${\kappa}$ then
$$
V[H]\models \text{``}\pcal({\omega})\cap V[H\rest I]\cap V[H\rest J]=
\pcal({\omega})\cap V.\text{''}
$$
\end{lemma}

\Proof{lemma}
Assume that $\adot$ is a nice $\coh I;$-name, $\bdot$ is a nice
$\coh J;$-name, $p\in \coh {\kappa};$ and
$p\force$ ``$\adot=\bdot$''. We can assume that $\dom(p)\subs I\cup J$.
We show that  for each $n\in {\omega}$ we have that 
$p\rest I\force$ ``$\hat{n}\in \adot$'' or
$p\rest I\force$ ``$\hat{n}\notin \adot$''. Indeed, if 
$p\rest I\notforce$ ``$\hat{n}\in \adot$'' then there is a condition
$q\le p\rest I$ in $\coh I;$ such that $q\force$ 
``$\hat{n}\notin \adot$'' and so $p\cup q\force$ 
``$\hat{n}\notin \bdot$''. But $\bdot$ is a $\coh J;$ name so
$(p\cup q)\rest J=p\rest J$ forces the same statement, 
$p\rest J\force$ 
``$\hat{n}\notin \bdot$''. 
But $\force$ $\adot=\bdot$ and so $p\force$ $ \hat{n}\notin \adot$ as well.
Thus $p$ decides the elements of $\adot$, in other words, 
$p\force$ ``{\em $\adot\in V$}''.
\Eproof


To conclude the proof we show that 
if  $q\in\coh {\kappa};$ then
there is a condition $r\le q$ in $\coh {\kappa};$ such that
$r\force $ ``\ref{pr:e}\ref{f1} or \ref{pr:e}\ref{f2} holds''.
Let $I'=I\cup \dom (q)$. 

For each $t\in T$ let ${\varphi}_t(x_0,\dots,x_{|t|-1})$ be the following
formula:
$$
{\varphi}(\<B_{0,k}:k<{\omega}\>,\dots,\<B_{|t|-1,k}:k<{\omega}\>)
\Longleftrightarrow \bigcap_{i<|t|} B_{i,t(i)}\in (\pcal({\omega}))^V.
$$
Applying the homogeneity lemma \ref{lm:homog} to $V[G\rest I']$ 
as our ground model we get that (A) or (B) below holds:
\begin{enumerate}\Alabel
\item 
$1_{\coh {\kappa};}\force$ ``
$\bint s;t;\in (\pcal({\omega}))^V$  for each $t\in T$ and 
$s\in\injs S;|t|;$,''  
\item 
for  some $t\in T$ and $p\in \coh {\kappa};$ we have \newline 
$p\force$ ``
there are subsets $\dot D_0$, $\dot D_1$, $\dots$, $\dot D_{|t|-1}$ of $S$
 such that
\begin{enumerate}\alabel
\item for each $A\in\br S;{\omega};\cap V$ we have  $A\cap\dot D_i\ne\empt$ 
for each $i<|t|$, 
\item If $s\in \injc \dot D_0, \dot D_1,\dots,  \dot D_{|t|-1};$  
 we have
$$
\dot\bint s;t;\notin (\pcal({\omega}))^V\text{.''}
$$
\end{enumerate}
\end{enumerate}

Let
$J_{\alpha}=\bigcup_{i<{\omega}}\supp(\bbbdot i;{\alpha};)$
for ${\alpha}\in S$ and 
 denote by 
$\fab$  the natural order preserving bijection between $J_{\alpha}$
and $J_{\beta}$ for $\<{\alpha},{\beta}\>\in  S^2$. 

Assume first that (A) holds.
Fix  $t\in T$ and
  $s\in \injs S;|t|;$. Write ${\alpha}_i=s(i)$ for $i<|t|$. Since
$\coh \kappa;$ is c.c.c , there is in $V$ a countable set 
$\ical_t\subs \pcal({\omega})$ such that
$1_{\coh {\kappa};}\force$ 
``$\bigcap\limits_{i<k}\bbbdot n_i;{\alpha}_i;\in \ical_t$''.

Assume that $\<{\delta}_0,\dots,{\delta}_{|t|-1}\>\in \injs S;k;$.

Let $J^*=\bigcup\limits_{i<k}J_{{\delta}_i}$, 
$J=\bigcup\limits_{i<k}J_{{\alpha}_i}$, and ${\psi}=\bigcup\limits_{i<k}{\varphi}_{{\delta}_i,{\alpha}_i}$.
Then ${\psi}$ is a bijection between $J^*$ and $J$ and so it lifts up
to an isomorphism between $\coh J^*;$ and $\coh J;$ and between the families
of nice $\coh J^*;$-names and  nice $\coh J;$-names.

Let $G$ be  $\coh {\kappa};$-generic and put $G_0=G\rest J^*$. 
Since $\supp(\bbbdot n_i;{\delta}_i;)\subs J^*$  it follows that 
$\ig G;{\bbbdot n_i;{\delta}_i;};=\ig G_0;{\bbbdot n_i;{\delta}_i;};$.
Let $G_1={\psi}''G_0$. Then $G_1$ is also 
a $\coh J;$-generic filter  and
since ${\psi}(\bbbdot n_i;{\delta}_i;)=\bbbdot n_i;{\alpha}_i;$ it follows that
\begin{equation}
\tag{$\bullet$}
\ig G_0;{\bbbdot n_i;{\delta}_i;};=\ig G_1;{\bbbdot n_i;{\alpha}_i;};.
\end{equation}

Since 
$1_{\coh {\kappa};}\force $ ``
$\bigcap\limits_{i<k}\bbbdot n_i;{\alpha}_i;\in \ical_t$'',  
by $(\bullet)$ we have  $1_{\coh {\kappa};}\force $
``$\bigcap\limits_{i<k}\ig G;{\bbbdot n_i;{\delta}_i;};\in \ical_t$''  
as well. From this it is obvious that we have
$$
1_{\coh {\kappa};}\force \{\dot\bint t;s;:t\in T\land
s\in \injs S;|t|;\}\subs \ical=\bigcup\{\ical_t:t\in T\}. 
$$
where $\ical$ is countable as $T$ is. 

Assume now that (A) fails and so (B) holds.

Let $G$  be $\coh {\kappa};$-generic with $p\in G$ and  
$\<{\gamma}_0,\dots,{\gamma}_{k-1}\>$, 
$\<{\delta}_0,\dots, {\delta}_{k-1}\>\in \injc D_0,\dots, D_{k-1};$
such that 
$$
V[G]\models \text{``$\{{\gamma}_i,{\delta}_i\}\in \br D_i;2;$ 
for $i<k$ are pairs  of distict ordinals''.}
$$

Let $J^*=\bigcup\limits_{i<k}J_{{\gamma}_i}$ and
 $J^{\star}=\bigcup\limits_{i<k}J_{{\delta}_i}$.
Then $J^*\cap J^\star=\empt$, hence by lemma \ref{lm:V} we have 
 $\pcal({\omega})\cap V[G\rest J^*]\cap\rest V[G\rest J^{\star}]=
\pcal({\omega})\cap V$ and so  
$V[G]\models $ ``$\bigcap_{i<k}\bbb n_i;{\delta}_i;\notin V$'' implies
that
$V[G]\models $ 
``$\bigcap_{i<k}\bbb n_i;{\delta}_i;\ne \bigcap_{i<k}\bbb n_i;{\gamma}_i;$''


The theorem is proved.
\eproof

\section{Applications}\label{sc:applications}

We start with  presenting some combinatorial applications because they are
quite simple and so they nicely illustrate the use of our principles.

Kunen \cite{Ku2} 
proved that if one adds  Cohen reals to a model of $CH$ then
in the generic extension there is no 
strictly $\subs^*$-increasing chain of subsets of ${\omega}$ of length $\oot$. 
The first theorem we prove easily yields a corollary which
is a generalization of Kunen's above result. 

\begin{theorem}\label{tm:tower}
If $\cccmin {\kappa};$ holds then for each 
$\acal\subs \br {\omega};{\omega};$ of size ${\kappa}$
either \begin{enumerate}\alabel
\item $\exists \bcal\in\br \acal;{\kappa};$ $\forall B\ne B'\in \bcal$
$|B\setm B'|={\omega}$  
\end{enumerate}
or
\begin{enumerate}\alabel \addtocounter{enumi}{1}
\item $\exists X\in \br {\omega};{\omega};$ 
$|\{A\in\acal:A\subs X\}|=|\{A\in\acal:X\subs^* A\}|={\kappa}$.
\end{enumerate}
\end{theorem}

\proof
Fix a 1--1 enumeration  $\{A_{\xi}:{\xi}<{\kappa}\}$ of $\acal$.
Let $A({\xi},2n)=A_{\xi}\setm n$ and 
$A({\xi},2n+1)=({\omega}\setm A_{\xi})\setm n$.
Put $T=\{\<2i,2i+1\>:i\in {\omega}\}$.
If $S\in \br {\kappa};{\kappa};$ witnesses  \ref{pr:c}\ref{c1}, then 
$\bcal=\{A_{\xi}:{\xi}\in S\}$ satisfies (a).
If on the other hand $D,E\in\br {\kappa};{\kappa}; $, $D\cap E=\empt$, with
  $\<2i,2i+1\>\in T$
show that \ref{pr:c}\ref{c2} holds, then let $X=\cup\{A_{\xi}:{\xi}\in D\}$.
Then $A_{\xi}\subs X$ for each ${\xi}\in D$ and
$X\setm i\subs A_{\zeta}$ for each ${\zeta}\in E$.
\eproof

\begin{corollary}
If $\cccmin \kappa;$ holds then there is no strictly 
$\subs^*$-increasing chain
 $\tcal\subs\br \omega;\omega;$
of length $\kappa$.
\end{corollary}

The next theorem can be considered as a kind of  dual to \ref{tm:tower}.

\begin{theorem}\label{tm:luzinplus}
If $\cccmin {\kappa};$ holds then for each 
$\acal\subs \br {\omega};{\omega};$ of size ${\kappa}$ and for each 
natural number $k$
either \begin{enumerate}\alabel
\item there is a family $\bcal\in\br \acal;{\kappa};$ such that 
for each $\bcal'\in \br \bcal;k;$ we have 
$|\bigcap\bcal'|={\omega}$  
\end{enumerate}
or
\begin{enumerate}\alabel \addtocounter{enumi}{1}
\item there are $k$ subfamilies $\bcal_0$, $\dots$, $\bcal_{k-1}$
of $\bcal$ of size ${\kappa}$ such that 
$$
|\bigcap_{i<k}\bigcup\bcal_i|<{\omega}.
$$
\end{enumerate}
\end{theorem}

\proof
Fix a 1--1 enumeration   $\{A_{\xi}:{\xi}<{\kappa}\}$ of $\acal$.
Let $A({\xi},n)=A_{\xi}\setm n$ and put
$T=\omega^k$.
If $C\in \br {\kappa};{\kappa};$ witnesses  \ref{pr:c}\ref{c1}, then 
$\bcal=\{A_{\xi}:{\xi}\in C\}$ satisfies (a).
If $D_0,\dots,D_{k-1}\in\br {\kappa};{\kappa}; $  and $\nok\in T$ 
show that \ref{pr:c}\ref{c2} holds, then we can assume that
the $D_i$ are pairwise disjoint and if we set 
$\bcal_i=\{A_\xi:\xi\in D_i\}$ then we have
$$
\bigcap_{i<k}\bigcup\bcal_i\subs \max\limits_{i<k}n_i.
$$
\eproof

\begin{remark}
In theorem \ref{tm:luzinplus}
we can not replace (a) with the following 
(slightly stronger) condition (a'):
\begin{enumerate}\alabel
\item[(a')] there is a family $ \bcal\in\br \acal;{\kappa};$ such that
for each 
$\bcal'\in \br \bcal;k+1;$ we have 
$|\bigcap\bcal'|={\omega}$,  
\end{enumerate}
and if $k>2$ then (b) can not be replaced by
\begin{enumerate}\alabel
\item[(b')] there are pairwise disjoint subsets 
$X_0,\dots,X_{k-1}$ of ${\omega}$ such that for each 
$i<k$ we have $\{A\in\acal:A\subs^* X_i\}\ne\empt$,
\end{enumerate}
because for each $k\in {\omega}$ one  can construct in ZFC a family 
$\acal\subs \br \omega;\omega;$
of size $2^{\omega}$ such that 
$\bigcap\acal'$ is finite for every 
$\acal'\in \br \acal;k+1;$
but $\bigcap\acal'$ is infinite whenever $\acal'\in \br \acal;k;$.  
Indeed, let $T=2^{<{\omega}}$ be the Cantor tree, 
and for $n<{\omega}$ let $C_n=2^n$
be the $n^{\rm th}$-level of $T$. 
For each $f\in 2^{\omega}$ let 
$$
A(f)=\bigcup_{n<{\omega}}\{X\in \br C_n;k;: f\rest n\in X\}
$$
and  $\acal=\{A(f):f\in 2^{\omega}\}$.
If $\bcal'\subs 2^{\omega}$ and $n<\omega$ then
$\bigcap\{A(f):f\in\bcal'\}\cap \br C_n;k;\ne\empt$ iff
 $|\{f\rest n:f\in\bcal'\}|\le k$. 
Thus $\acal$ satisfies our requirements.
This example is due to A. Hajnal and included here with his kind
permission.
\end{remark}

Next we prove  a consequence of theorem \ref{tm:luzinplus},
but first we give a definition.

\begin{definition}
Let ${\kappa} $ be a regular cardinal and 
$\acal \subs\br {\omega}; {\omega} ; $ be an almost disjoint
family.  $\acal$ is called a {\em $ {\kappa}  $-Luzin gap } 
if $|\acal|={\kappa}$ and there is no $X\in\br {\omega} ;{\omega} ;$ 
such that both $|\{A\in\acal:|A\setm X|<{\omega}|\}={\kappa} $ and   
$|\{A\in\acal:|A\cap X|<{\omega}\}|={\kappa} $.   
A {\em Luzin-gap} is an ${\omega_1} $-Luzin gap.
\end{definition}

An $\oo$-Luzin gap can be  constructed in ZFC and simple
forcings give models in which  there are $2^ {\omega} $-Luzin gaps
while $2^{\omega}$ is as large as you wish.
The next corollary of theorem \ref{tm:luzinplus} implies that 
one can not construct  ${\omega}_2$-Luzin gaps  from the assumption
$2^{\omega}\ge{\omega}_2$ alone .

\begin{corollary}\label{cor:luzingap}
If $\cccmin {\kappa};$ holds then there is no ${\kappa}$-Luzin gap.
\end{corollary}

\proof
Assume that 
$\acal\subs \br {\omega};{\omega};$ is an almost disjoint family of
 size ${\kappa}$. 
Then we can not get a even a two element subfamily 
$\bcal\subs\acal$  satisfying \ref{tm:luzinplus}(a). 
So applying theorem \ref{tm:luzinplus} for 
this $\acal$ and for $k=2$   there are subfamilies 
$\bcal\subs\acal$ and $\dcal\subs\acal$ of size ${\kappa}$ 
such that $(\bigcup\bcal)\cap (\bigcup\dcal)$ is finite.
Hence $X=\bigcup\bcal$ witnesses that $\acal$ is not a 
${\kappa}$-Luzin gap.
\eproof

We have one more theorem of this type.

\begin{theorem}\label{tm:luzinplus-omega}
If $\cccmin {\kappa};$ holds then for each 
$\acal\subs \br {\omega};{\omega};$ of size ${\kappa}$ 
either \begin{enumerate}\alabel
\item there is a centered subfamily $\bcal\subs\acal$  
of size ${\kappa}$,
\end{enumerate}
or
\begin{enumerate}\alabel \addtocounter{enumi}{1}
\item for some $k<{\omega}$ 
there are subfamilies $\bcal_0$, $\dots$, $\bcal_{k-1}$
of $\bcal$ of size ${\kappa}$ such that 
$$
|\bigcap_{i<k}\bigcup\bcal_i|<{\omega}.
$$
\end{enumerate}
\end{theorem}

\proof
We can argue as in the proof of  theorem \ref{tm:luzinplus}
using $T={\omega}^{<{\omega}}$ instead of $T={\omega}^k$.
\eproof

For $\acal\subs \pcal({\omega})$ and  ${\kappa}<{\omega}$ let 
$\pint \acal;k;=\{\cap{\acal}':{\acal}'\in\br \acal;n;\}\}$. 
Put $\pint \acal;<{\omega};=\bigcup\limits_{k<{\omega}}\pint \acal;k;$.

\begin{theorem}\label{tm:int-sizek}
If  $\eeemin {\kappa};$ holds then  for each family   
$\acal\subs \pcal({\omega})$ 
of size ${\kappa}$ and for each natural number $k$  
either 
\begin{enumerate}\alabel
\item $|\pint \acal;k;|={\kappa}$
\end{enumerate}
or
\begin{enumerate}\alabel \addtocounter{enumi}{1}
\item there is a subfamily $\bcal\subs \acal$ of size ${\kappa}$ such that
$|\pint \bcal;k;|\le {\omega}$.
\end{enumerate}
 \end{theorem}

\proof
Fix a 1--1 enumeration $\{A_{\alpha}:{\alpha}<{\kappa}\}$ of $\acal$,
let $T={\omega}^k$
and consider the matrix 
$\acal'=\<\aaa n;{\alpha};:{\alpha}<{\kappa};n<{\omega}\>\in \matr {\kappa};$ 
defined by the stipulation $\aaa n;{\alpha};=A_{\alpha}$.
Apply $\eeemin {\kappa};$. If \ref{pr:e}(2) holds, then 
$|\pint \acal;k;|={\kappa}$. If $S\in \br {\kappa};{\kappa};$ 
witnesses \ref{pr:e}(1) then subfamily 
$\bcal=\{A_{\alpha}:{\alpha}\in S\}$ satisfies 
$|\pint \bcal;k;|={\omega}$.
\eproof

\begin{theorem}\label{tm:int-sizeo}
If  $\eeemin {\kappa};$ holds then  for each family   
$\acal\subs \pcal({\omega})$ 
of size ${\kappa}$   
either 
\begin{enumerate}\alabel
\item there is a natural number $k$ such that $|\pint \acal;k;|={\kappa}$
\end{enumerate}
or
\begin{enumerate}\alabel \addtocounter{enumi}{1}
\item there is a subfamily $\bcal\subs \acal$ of size ${\kappa}$ such that
$|\pint \bcal;<{\omega};|\le {\omega}$.
\end{enumerate}
 \end{theorem}

\proof
We can argue as in the proof of  theorem \ref{tm:int-sizek}
using $T={\omega}^{<{\omega}}$ instead of $T={\omega}^k$.
\eproof

Now we turn to applying our principles to topology.
We start with an application of the relatively weak principle
$\diamin {\kappa};$.

A. Dow \cite{Dow2} proved that if we add ${\omega}_2$ Cohen
reals to a model of $GCH$ then in the generic extension
$ {\beta}{\omega} $  can be embedded into 
every separable, compact $T_2$ space of size
$>c={\omega}_2$. Here we show that $c=\oot=2^{\oo}$ together with  
$\diamin {\omega}_2;$  suffice to imply this statement.

First we need a lemma based on the observation that  large separable spaces 
contain many ``similar'' points. 

Given a topological space $X$ and a point $x\in X$ we denote by 
$\vcal_X(x)$ the neighbourhood filter of $x$ in $X$, that is,
$\vcal_X(x)=\{U\subs X:x\in \operatorname{int}_X(U)\}$.
If $D$ is a dense subset of $X$ let 
$\vcal_X(x)\rest D=\{U\cap D:D\in\vcal_X(x)\}$.
We  omit the subscript $X$ if it may not cause any confusion.

In section \ref{sc:principles} we defined the operation
$\lsup \acal;{\omega};$ 
for $\acal\in \matr {\kappa};$. 
By an abuse of notation 
we define $\lsup \acal;{\omega};$
 for every family $\acal$ of subsets of 
${\omega}$  as follow:
$$
\lsup\acal;{\omega};=\{X\subs {\omega}:|\acal\cap P(X)|=|\acal|\}.
$$

\begin{lemma}\label{lm:bigx}
Assume that $X$ is a separable regular topological space of size $>c^{<c}$,
where $c=2^{\omega}$, $D\in \br X;{\omega};$, $\overline{D}=X$. 
Then there are  a point
$x\in X$ and  a family 
$\acal=\{A_{\alpha}, B_{\alpha}:{\alpha}<c\}\subs \pcal(D)$ such that
\begin{enumerate}\arablabel
\item $\overline{A_{\alpha}}\cap \overline{B_{\alpha}}=\empt$ for each 
$\alpha<c$,
\item $\lsup \acal;D;\subs \vcal(x)\rest D$.
\end{enumerate}
\end{lemma}

\proof
Fix an enumeration $\{D_{\xi}:{\xi}<c\}$ of $\pcal(D)$ and
let $\dcal_{\alpha}=\{D_{\xi}:{\xi}<{\alpha}\}$
 for ${\alpha}<c$.
For $x\in X$ and ${\alpha}<c$ let
$\vcal(x,{\alpha})=(\vcal(x)\rest D)\cap \dcal_{\alpha}$.
A point $x\in X$ is called {\em special} if there is an ${\alpha}<c$
such that $\vcal(x,{\alpha})\ne\vcal(y,{\alpha})$ for each 
$y\in X\setm\{x\}$. Clearly there are at most $c^{<c}$ special
points in $X$. Since $|X|>c^{<c}$ we can pick a point 
$x\in X$ which is not special. Then for each ${\alpha}<c$ we can
find a point $x_{\alpha}\ne x$ in $X$ such that 
$\vcal(x_{\alpha},{\alpha})=\vcal(x,{\alpha})$. Since $X$ is regular the points
$x$ and $x_{\alpha}$ have neighbourhoods $U_{\alpha}$ and $W_{\alpha}$,
respectively, with $ \overline{U_{\alpha}}\cap \overline{W_{\alpha}}=\empt$.
Let $A_{\alpha}=U_{\alpha}\cap D$ and $B_{{\alpha}}=W_{\alpha}\cap D$.

Now assume that $E\in\lsup \acal;D;$ and pick   ${\xi}<c$ with
$E=D_{\xi}$. We can find ${\alpha}<c$ such that ${\xi}<{\alpha}$ and
either $A_{\alpha}\subs E$ or $B_{\alpha}\subs E$. Hence
$E\in\vcal(x,{\alpha})\cup\vcal(x_{\alpha},{\alpha})=\vcal(x,{\alpha})$.
Therefore $E\in\vcal(x)\rest D$ which was to be proved.
\eproof

Let us now recall the definition of  a $ {\mu}$-dyadic system
from \cite{J}. 
\begin{definition}
If $X$ is a topological space
a  family 
$\{\<\aaa 0;{\alpha};,\aaa 1; {\alpha};\> :{\alpha}\in {\mu}\}$ of
pairs of closed subsets of $X$  is a
 {\em ${\mu}$-dyadic system}   such that
\begin{enumerate}
\item $\aaa 0;{\alpha};\cap \aaa 1; {\alpha};=\empt$ 
for each ${\alpha}<{\mu}$,
\item  for each ${\epsilon}\in \fn({\mu},2,{\omega})$
we have $\bigcap\limits_{{\alpha}\in\dom({\epsilon})} 
\aaa {\epsilon}({\alpha}); {\alpha};\ne\empt$. 
\end{enumerate}
\end{definition}

\begin{theorem}\label{tm:betaomega}
If  $\diamin c;$ holds, $X$ is a  separable compact  
$T_2$ space of size
$>c^{<c}$ then $X$ contains a $c$-dyadic system, consequently 
$X$ maps continuously onto $[0,1]^c$ $($
and so $\beta\omega$ can be embedded  into $X$ $)$.
\end{theorem}

\proof
Fix a countable dense subset $D$ of $X$.
By lemma \ref{lm:bigx} there is a family 
$\acal=\{A_{\alpha}, B_{\alpha}:{\alpha}<c\}\subs \pcal(D)$ 
such that
$\overline{A_{\alpha}}\cap \overline{B_{\alpha}}=\empt$ for $\alpha<c$ and
 $\lsup \acal;D;$ is centered. 
Let $\ddd 0;{\alpha};=A_{\alpha}$, $\ddd 1;{\alpha};=A_{\alpha}$ and
$\ddd n;{\alpha};=D$ for ${\alpha}<{\kappa}$ and $n\ge 2$
and consider the ${\kappa}\times {\omega}$-matrix 
$\dcal=\<\ddd i;{\alpha};:{\alpha}<{\kappa},i<{\omega}\>$. 
Since $\lsup \acal;D;=\lsup\dcal;D;$ we can apply $\diamin c;$ to get 
a cofinal $S\subs c$ such that the family   
$\<\overline{A_{\alpha}}, \overline{B_{\alpha}}:{\alpha}<c\>$
is $c$-dyadic. Now we can apply theorem \cite[3.18]{J} 
to get the other consequences. 
\eproof

A topological space $X$ is called {\em scattered} if every subspace of $X$
has an isolated point. 
For a scattered  space  $X$ 
denote by $\levx {\alpha};$ the ${\alpha}^{\text{th}}$ 
Cantor-Bendixon level of $X$.
The height of $X$, $\htt(X)$,
is defined as the minimal ${\alpha}$ with $\levx {\alpha};=\empt$.
Following  \cite{Ro} we call $X$  {\em thin} if all  levels of $X$
are  countable.

Since the cardinality of a locally compact, scattered separable space
is at most $2^{\omega}$ by \cite{MRS}, the height of such a space
is less then $({2^{\omega}})^{{}^+}$. So under $CH$ there is no such a space
of height ${\omega}_2$.  
I. Juh\'asz and W.  Weiss, \cite[theorem 4]{JW}, proved in ZFC that
for every ${\alpha}<{\omega}_2$ there is a locally compact, scattered
thin space $X$ with $\htt(X)={\alpha}$. 
M. Weese asked  whether  the existence of 
such a space of height ${\omega}_2$ follows from $\neg CH$.
This question was answered in the negative  by W. Just, who proved,
\cite[theorem 2.13 ]{J1}, that if one adds Cohen reals 
to a model of $CH$ then in the generic extension there are no
locally compact scattered thin spaces of height ${\omega}_2$.
On the other hand,
J. Baumgartner and S. Shelah, \cite{BS}, constructed a ZFC model which
contains such a space of height ${\omega}_2$. 

The next theorem is a generalization of the above mentioned
result of Just.

\begin{theorem}\label{tm:thin-tall}
If $\ccc {\kappa};$ holds then there is no locally compact, thin
scattered space of height ${\kappa}$. 
\end{theorem}

\def\uu#1;#2;{U(#1,#2)}
\def\bb#1;#2;{B(#1,#2)}
\def\cc#1;#2;{C(#1,#2)}
\def\aa#1;#2;{A(#1,#2)}
\def\dddd#1;#2;{D(#1,#2)}

\proof
Assume on the contrary that there is such a space $X$. We can assume that
 $\levx {\alpha};=\{{\alpha}\}\times {\omega}$ for ${\alpha}<\htt(X)$.
For each ${\alpha}<\htt(X)$ fix  compact open  neighbourhoods 
$\uu {\alpha};n;$ of $\<{\alpha},n\>$ for $n\in {\omega}$ such that
$\uu \alpha;n;\subs \{\<\alpha,n\>\}\cup\bigcup\{\levx \beta;:\beta<\alpha\}$
and 
the sets $\uu {\alpha};n;$ for $n<{\omega}$ are pairwise disjoint.

Put $\aa {\alpha};2n;=\uu {\alpha};n;\cap \levx 0;$ and 
$\aa {\alpha};2n+1;=\levx 0;\setm\bigcup\{\uu {\alpha};m;:m\le n\}$.
Let 
$$T=\{t\in {\omega}^{<{\omega}}:\text{
$t(0)$ is even and $t(i)$ is odd for $i>0$
}\}.
$$

Now apply $\ccc {\kappa};$ to the matrix 
$\<\aaa n;{\alpha};:{\alpha}<{\kappa},n<{\omega}\>\in\matr {\kappa};$
and $T$.

\begin{tabular}{ccl}
Observe that 
&$\aa {\beta};2n;\cap\bigcap\limits_{i<k}\aa {\alpha}_i;2n_i+1;=\empt$&iff\\
&$\uu {\beta};n;\cap \levx 0;
\subs\bigcup\limits_{i<k}\uu {\alpha}_i;n_i;\cap \levx 0;$&iff\\
&$\uu {\beta};n;\subs\bigcup\limits_{i<k}\uu {\alpha}_i;n_i;$.&
\end{tabular}

Thus if  $t=\<2n,2n_0+1,\dots,2n_{k-1}+1\>\in T$ and 
$\<{\beta},{\alpha}_0,\dots, {\alpha}_{k-1}\>\in  \injs {\kappa};k+1;$ then
$\aa {\beta};2n;\cap \bigcap\limits_{i<k}\aa {\alpha}_i;2n_i+1;=\empt$
implies ${\beta}\le\max\limits_{i<k}{\alpha}_i$. This excludes 
\ref{pr:c}(2).
So  \ref{pr:c}(1) holds, that is we have a stationary set
$S\subs {\kappa}$ such that if $t=\<2n,2n_0+1,\dots,2n_{k-1}+1\>\in T$
and $\<{\beta},{\alpha}_0,\dots, {\alpha}_{k-1}\>\in \injs S;k+1;$ then
$$
\aa {\beta};n;\cap \bigcap\limits_{i<k}\aa {\alpha}_i;2n_i+1;\ne\empt,
$$
that is 
$$
\uu {\beta};n;\setm\bigcup\limits_{i<k}
\bigcup\limits_{j\le n_i}\uu {\alpha}_i;j;\ne\empt.
$$
But $\uu {\beta};n;$ is compact and each $\uu {\alpha};n;$ is open 
so it follows that for every ${\beta}\in S$ and $n\in {\omega}$  the set
$$
\dddd {\beta};n;=\uu {\beta};n;\setm
\bigcup\{\uu({\alpha};m;:{\alpha}\in S\setm\{{\beta}\}\land m\in {\omega}\}
$$
is not empty.
For every such ${\beta}$ and $n$
let $\<{\gamma}({\beta},n),m({\beta},n)\>\in \dddd {\beta};n;$.

Since $\levx {\beta};$ is dense in 
$X\setm \bigcup\{ \levx {\alpha};:{\alpha}<{\beta}\}$ for every
${\beta}\in {\kappa}$
there is $k({\beta})\in {\omega}$ such that 
$\<{\beta},k({\beta})\>\in \uu {\beta}^*;0;$, 
where ${\beta}^*=\min S\setm {\beta}+1$.
Thus $\<{\beta},k({\beta})\>\notin \dddd {\beta};k({\beta});$ and so
${\gamma}({\beta},k({\beta}))<{\beta}$ for each  ${\beta}\in S$.
The set $S$ is stationary so there are  a stationary set $S'\subs S$, 
and ordinals ${\gamma}<{\kappa}$ and $k,m<{\omega}$ such that 
$k({\beta})=k$,
${\gamma}({\beta},k)={\gamma}$ and $m({\beta},k)=m$ whenever ${\beta}\in S'$.
Thus $\<{\gamma},m\>\in \dddd {\beta};k;$ for each ${\beta}\in S'$,
while $\dddd {\beta};k;\cap \dddd {{\beta}'};{k};=\empt$ for any
$\{\beta,\beta'\}\in\br S';2;$ by the construction.
This is a contradiction, hence the theorem is proved.
\eproof

In \cite{J1} W. Just also proved that if one adds at least
${\omega}_2$ Cohen reals to a model of $CH$ then in the generic extension
there is no locally compact, scattered
topological space $X$ such that $\htt(X)={\oo}+1$,
$\levx 0;$ is countable, $|\levx {\alpha};|\le{\oo}$ for ${\alpha}<{\oo}$
and $|\levx {\oo};|={\omega}_2$.
The next theorem shows how to  get a generalization of this result
from our principles.

\begin{theorem}
If $\cf({\lambda})\ge\oo$ and  $\eeemin {\lambda}^+;$ holds then 
there is no locally compact, scattered
topological space $X$ such that $\htt(X)={\lambda}+1$,
$\levx 0;$ is countable, $|\levx {\alpha};|\le{\lambda}$ for all 
${\alpha}<{\lambda}$ and $|\levx {\lambda};|={\lambda}^+$.
\end{theorem}

\proof
Assume on the contrary that $X$ is such a space.

We can assume that $\levx 0;={\omega}$ and that 
$\levx {\lambda};=\{{\lambda}\}\times {\lambda}^+$.
For each $x\in X$ choose a compact open neighbourhood
$U(x)$ of $X$ and let $B(x)=U(x)\cap {\omega}$. 
Put $\bcal=\{B(x):x\in \levx {\lambda};\}$.
Let 
$\ucal=\operatorname{CO}(X\setm \levx {\lambda};)$, i.e. the
family of compact open subsets of 
$X\setm \levx {\lambda};=\bigcup\{\levx \alpha;:\alpha<{\lambda}\}$. 
and $\fcal=\{U\cap {\omega}:U\in \ucal\}$.
Since $X$ is locally compact it follows that for each 
$\{x,y\}\in\br {\levx {\lambda};};2;$ we have
$U(x)\cap U(y)\in \ucal$ and so $B(x)\cap B(y)\in\fcal$.

Since $|\fcal|\le {\lambda}$, it follows 
$|\pint \bcal;2;|\le {\lambda}<{\lambda}^+$.
Thus, applying theorem \ref{tm:int-sizek} for $k=2$ we can get 
 a cofinal set $S\subs {\lambda}^+$ such that
the family
$$
\ical=\{B(\<{\lambda},{\alpha}\>)\cap B(\<{\lambda},{\beta}\>):
\{{\alpha},{\beta}\}\in \br S;2;\}
$$
is at most countable.

Then there is ${\gamma}<{\lambda}$ such that $\ical$ is  contained in 
$\operatorname{CO}(\bigcup\{\levx {\gamma}';:{\gamma}'<{\gamma}\})$.
Therefore $U(y)\cap U(y')\subs\bigcup\{\levx {\gamma}';:{\gamma}'<{\gamma}\}$
for each $\{y,y'\}\in \br \{{\lambda}\}\times S;2;$
and so the sets $U(x)\cap \levx {\gamma};$  for 
$x\in \{{\lambda}\}\times S$ are pairwise disjoint and non-empty
which contradicts $|\levx {\gamma};|\le {\lambda}$. 
\eproof

Following the terminology of \cite{HJ}
a Hausdorff space is called $P_2$ if   
it does not contain two uncountable disjoint open sets.
Hajnal and Juh\'asz in \cite{HJ} constructed a ZFC example of a 
first countable, $P_2$ space of size $\oo$
as well as consistent examples of size $2^{\omega}$ with 
$2^{\omega}$ as large as you wish . On the other hand, using a result
of Z. Szentmikl{\'o}ssy they proved  that it is consistent with ZFC
 that $2^{\omega}$ is as large as you wish and
there are no first countable $P_2$ spaces of size 
$\ge {\omega} _3$. However their method was unable to 
replace here ${\omega}_3$ with ${\omega}_2$.
Our next result does just this because, as is shown in \cite{HJ},
every $P_2$ space is separable.
 
\begin{theorem}
\label{tm:p2}
If $\cccmin {\kappa};$ holds then every
 first countable, separable $T_2$ topological space $X$ of size ${\kappa}$
contains two disjoint open sets $U$ and $V$ of cardinality ${\kappa}$.
\end{theorem}

\proof
Let $D$ be a countable dense subset of $X$. For each $x\in X$
fix a neighborhood base $\{U(x,n):n\in {\omega}\}$ of $x$ in $X$.
Apply $\cccmin{\kappa};$ to the matrix 
 $\<U(x,n)\cap D:x\in X, n<{\omega}\>$ and  $T={}^2{\omega}$. 
Since $X$ is $T_2$, there is no 
$S\in\br X;{\kappa};$ satisfying \ref{pr:c}\ref{c1}. So there are 
$S_0,S_1\in \br X;{\kappa};$ and $n,m\in {\omega}$ such that 
$U(x,n)\cap U(y,m)\cap D=\empt$ whenever $x\in S_0$ and $y\in S_1$.
But $D$ is dense in $X$, therefore $U=\bigcup\{U(x,n):x\in S_0\}$ and
$V=\bigcup\{U(y,m):y\in S_1\}$ are disjoint open sets of size ${\kappa}$.
\eproof

\begin{definition}\label{df:seq_dense}
Let $X$ be a topological space and $D\subs X$. We say that
$D$ is {\em sequentially dense in $X$} iff for each $x\in X$ there is a 
sequence $S_x$ from $D$ which converges to $x$.
A space $Y$ is said to be {\em sequentially separable} if it contains a countable
sequentially dense subset.
\end{definition}

\begin{definition}
Given a topological space $\<X,{\tau}\>$ and a subspace $Y\subs X$ a function $f$
is called a {\em \nea on $Y$ in $X$} iff 
$f:Y\to {\tau}$ and $y\in f(y)$ for each 
$y\in Y$. 
\end{definition}

Our next result says that under $\cccmin {\kappa};$ if a sequentially 
separable space
$X$ does not contain a discrete subspace of size ${\kappa}$,
(i.e. $\shat(X)\le{\kappa}$ using the notation of 
\cite{J}) then $X$ does not contain  left or  
right separated subspaces of size ${\kappa}$ either. 
This   can be written as $\hhat(X)\zhat(X)\le{\kappa} $. 
Since in \cite{JS}  a  normal, Frechet-Urysohn, separable 
(hence sequentially separable)   space $X$  is forced such that   
$\operatorname{z}(X)\le \oo$ but $\hl(X)={\omega}_2$,  this
result  is not provable in ZFC. First, however, we need a lemma.

\begin{lemma}\label{lm:seq_sep}
Assume that $\cccmin {\kappa};$ holds.
Let $X$ be a sequentially separable space with  $Y\subs X$, $|Y|={\kappa}$.
If $f$ is a \nea on $Y$ in $X$, then either (a) or (b) below holds:
\begin{enumerate}\alabel
\item there is $Y'\in \br Y;{\kappa};$ such that 
$f(y)\cap Y'=\{y\}$ for each $y\in Y'$ (hence
$Y'$ is discrete),
\item there are $Y_0$, $Y_1\in \br Y;{\kappa};$ such that 
$y\in \overline{f(x)}$ whenever $x\in Y_0$ and $y\in Y_1$.
\end{enumerate}
\end{lemma}

\proof
We can assume that $D={\omega}$ is sequentially dense in $X$. For each
$y\in Y$ choose a sequence $S_y\subs D$ converging to $y$. 
Let $A(y,2n)=D\setm f(y)$,  $A(y,2n+1)=S_y\setm n$,
$T=\{\<2n,2m+1\>:n,m\in {\omega}\}$ and apply 
$\cccmin{\kappa};$. Assume first that   $Y'\in \br Y;{\kappa};$ 
witnesses \ref{pr:c}\ref{c1} and let $x\ne y\in Y'$. 
Then for each $n\in {\omega}$
we have $(S_y\setm f(x))\setm n\ne\empt$, i.e. $S_y\setm f(x)$ is infinite.
But $S_y$ converges to $y$, so $y\notin f(x)$, and so $Y'$ satisfies
(a). Assume now that \ref{pr:c}\ref{c2} holds. Then there are 
$Y_0,Y_1\in \br {\omega};{\omega};$
and $m\in {\omega}$ such that $(D\setm f(x))\cap (S_y\setm m)=\empt$ for each 
$x\in Y_0$ and $y\in Y_1$. But then $S_y\setm m\subs f(x)$ hence
$y\in\overline{f(x)}$ which was to be proved. 
\eproof

\begin{theorem}\label{tm:shz}
If $\cccmin{\kappa};$ holds, $X$ is a regular, sequentially separable space
with $\shat(X)\le {\kappa}$ then $\hhat(X)\zhat(X)\le{\kappa}$.
\end{theorem}

\proof
Assume on the contrary that $Y\in \br X;{\kappa};$ and the \nea 
$f:Y\to{\tau}_X$ witnesses that  $Y$ is left (right) separated. 
We
can assume that $Y={\kappa}$ and $Y$ is left (right) separated under the 
natural ordering of ${\kappa}$.
Since $X$ is regular we can find a \nea $g:Y\to{\tau}$ with
$\overline{g(y)}\subs f(y)$ for each $y\in Y$.
Apply  lemma \ref{lm:seq_sep} to $Y$ and $g$. 
Now \ref{lm:seq_sep}(a) can not hold because 
 $\shat(X)\le {\kappa}$, hence 
 there are $Y_0, Y_1\in \br Y;{\kappa};$ satisfying \ref{lm:seq_sep}(b).
Since both $Y_0$ and $Y_1$ are cofinal in $Y={\kappa}$ under the natural
ordering of the ordinals, applying left (or right) separatedness of $Y$
we can pick $x\in Y_0$ and $y\in Y_1$ such that $y\notin f(x)$.
By the choice of $g$  this implies $y\notin \overline{g(x)}$ which 
contradicts \ref{lm:seq_sep}(b).
\eproof

The Sorgenfrey line $L$  is weakly separated and is of size $c$ 
with $\shat(L)=\oo$. This shows 
 that  theorem \ref{tm:shz} does not remain valid if you put   
weakly separated subspaces instead of right or left separated ones.

As an easy consequence of \ref{tm:shz} we get the following result in which
(sequential) separability is no longer assumed.
We also note that under CH the assumption of
$X$ being Frechet-Urysohn is not necessary in this result.

\begin{theorem}
Assume $\cccmin \oot;$. If $X$ is regular, 
Frechet-Urysohn space and $\sprd(X)={\omega}$ then
$\hl(X)\le \oo$. 
\end{theorem}

\proof
If $\cccmin \oot;$ and
$X$ is separable, then by theorem  \ref{tm:shz} even  
$\sprd(X)\le \oo$ implies $\hl(X)\operatorname{z}(X)\le\oo$.
Now, every uncountable space $X$ which is both right and left 
separated contains an uncountable discrete subspace, hence 
every right separated subspace of $X$ is (hereditarily)
separable. So by the above if $Y\subs X$ is right separated then 
$|Y|\le \oo$, i.e. $\hl(X)\le \oo$.
\eproof

In  \cite{JSS2} we investigated the  following question:
What makes a space have weight larger than its character?
To answer this question we introduced the notion of  an
{\em irreducible base of a space}  and proved that
any  weakly separated space has such a base, moreover 
the weight of a space possessing an irreducible base 
can not be smaller than its cardinality.
We asked \cite[Problem 1]{JSS2} whether
 every  first countable space of uncountable weight 
contains an uncountable subspace with an irreducible base?
In theorem  \ref{tm:irred} and corollary 
\ref{cor:irred} we will give  a partial positive answer to this
problem, using the principle $\cccmin {\kappa};$.
First we recall some definitions from \cite{JSS2}.

\begin{definition}\label{df:irred}
Let $X$ be a topological space. 
A base $\ucal$ of $X$ is called {\em irreducible} 
if it has an {\em irreducible decomposition} $\ucal=\bigcup\{\ucal_x:x\in X\}$,
i.e.,  $(i)$ and $(ii)$ below hold:
\begin{enumerate}\rlabel
\item $\ucal_x$ is a neighbourhood base of $x$ in $X$ for each $x\in X$,
\item for each $x\in X$   the family 
$\ucal^-_x=\bigcup\limits_{y\ne x}\ucal_y$ is not a base of $X$
(hence  $\ucal^-_x$ does not contain a neighbourhood base of 
$x$ in $X$). 
\end{enumerate}
\end{definition}

\begin{definition}\label{df:irred-outer}
Let $X$ be a topological space with $Y\subs X$. 
Similarly as above, an outer base $\ucal$ of $Y$ in 
$X$ is called {\em irreducible} 
if it has an {\em irreducible decomposition} $\ucal=\bigcup\{\ucal_y:y\in Y\}$,
i.e.,  $(i)$ and $(ii)$ below hold:
\begin{enumerate}\rlabel
\item $\ucal_y$ is a neighbourhood base of $y$ in $X$ for each $y\in Y$,
\item for each $y\in Y$   the family 
$\ucal^-_y=\bigcup\{\ucal_z:z\in Y\setm \{y\}\}$ does not 
contain a neighbourhood base of $y$ in $X$. 
\end{enumerate}
\end{definition}

Note that in general, a subspace $Y$ having an irreducible outer base in $X$
does not necessarily possess an irreducible base in itself. However, if $Y$
is dense in an open set and the irreducible 
outer base of $Y$ consists of regular open sets then clearly this is the case.
Moreover, by our next result, under certain conditions we can at least find 
another subspace of the same size as $Y$ that does have an irreducible base.

\begin{lemma}\label{lm:regirr}
If $X$ is a regular, separable space and
 $Y\subs X$ has an irreducible outer base in $X$ 
consisting of regular open sets, then there is $Z\subs X$ with $|Z|=|Y|$
such that the subspace $Z$ has an irreducible base.
\end{lemma}

\proof
Let $\bcal=\bigcup\{\bcal_y:y\in Y\}$ be an irreducible outer base of $Y$ 
in $X$ consisting of regular open sets 
and  $D$ be a countable dense subset of $X$. We distinguish two cases:
\newcases
\begin{case}
$|(\operatorname{int}\overline{Y})\cap Y|=|Y|$.
\end{case} 
Let  $Z=(\operatorname{int}\overline{Y})\cap Y$. Since
$Z$ is dense in the open set $\operatorname{int}\overline{Y}$,  
by our above remark
$Z$  has an irreducible base.
\begin{case}
$|(\operatorname{int}\overline{Y})\cap Y|<|Y|$. 
\end{case} 
\noindent
In this case the set $Y_1=Y\setm \operatorname{int}\overline{Y}$ is 
nowhere dense,
so $D_1=D\setm \overline{Y_1}$ is dense in $X$. Let $Z=D_1\cup Y_1$, then
$|Z|=|Y_1|=|Y|$.
Write $D_1=\{d_n:n<\omega\}$ and for each $d_n\in D$ let 
$\bcal_{d_n}$ be a neighbourhood base of $d_n$ in $X$ consisting of 
regular open sets, that are  disjoint to 
$Y_1\cup\{d_m:m<n\}$. 
Then clearly 
$\bigcup\{\bcal_z:z\in Z\}$ is an irreducible outer base of $Z$ in $X$ consisting of
regular open sets and $Z$ is dense in $X$, so again we are done.
\eproof

\begin{theorem}\label{tm:irred}
Let ${\kappa}$ be a regular cardinal and assume $\cccmin{\kappa};$. 
If $X$ is a separable, first countable, 
regular space with $\w(X)\ge {\kappa}$, then there is subspace 
$Y\subs X$ of cardinality ${\kappa}$ that has an
irreducible  base. 
\end{theorem}

\proof
Let $D\subs X$ be a countable, dense subset of $X$. 
For each $x\in X$ fix a  neighbourhood 
base $\{U(x,n):n\in {\omega}\}$  consisting of regular open sets
and set  $V(x,n)=U(x,n)\cap D$. 
Since the $U(x,n)$ are regular open and $D$ is dense, we clearly 
have $U(x,n)\subs U(y,m)$ iff $V(x,n)\subs V(y,m)$.

Since $\w(X)\ge {\kappa}$, by transfinite recursion on $\beta<\kappa$
we can choose points $\{x_{\alpha}:{\alpha}<{\kappa}\}\subs X$
such that for any ${\beta}<{\kappa}$ the family $\{U(x_\alpha,n):\alpha<\beta,n<\omega\}$ does not contain
 a neighbourhood base of $x_\beta$, in other words, there is a
natural number $k_{\beta}$
such that for all ${\alpha}<{\beta}<{\kappa}$ and $n\in {\omega}$
we have 
\begin{equation}\tag{$*$}
\neg (x_{\beta}\in U(x_{\alpha},n)\subs U(x_{\beta}, k_{\beta})).
\end{equation}
We can assume that $k_{\beta}=0$ for each ${\beta}<{\kappa}$.
Let $X'=\{x_{\alpha}:{\alpha}<{\kappa}\}$. For $x\in X'$ and $n<{\omega}$
put 
$$
A(x,2n)=\left[V(x,n)\times\{0\}\right]\cup 
\left[(D\setm V(x,n))\times \{1\}\right]$$
and 
$$
A(x,2m+1)=\left[(D\setm V(x,0))\times\{0\}\right]\cup 
\left[V(x,m)\times \{1\}\right].
$$
Note that $A(x,2n)\cap A(y,2m+1)=\empt$ iff $V(y,m)\subs V(x,n)\subs V(y,0)$.
Apply  $\cccmin{\kappa};$ to $\<A(x,i):x\in X',i<\omega\>$ and
 $T=\{\<2n,2m+1\>:n,m<\omega\}$.
By $(*)$ (and $k_{\beta}=0$) there are no $D,E\in \br X';{\kappa};$ and $n,m\in {\omega}$ such that 
\begin{equation}\tag{$\dag$}
V(y,m)\subs V(x,n)\subs V(y,0)
\end{equation}
whenever $x\in D$ and $y\in E$, because $(\dag)$ fails if $x=x_{\alpha}$,
$y=x_{\beta}$ and ${\alpha}<{\beta}$.
So there is $Y\in \br X';{\kappa};$ such that for all $n,m\in {\omega}$
and $x\ne y\in Y$ the intersection of $A(x,2n)$ and $A(y,2m+1)$ is not empty.
This means that  $\neg (V(y,m)\subs V(x,n)\subs V(y,0))$, i.e. 
if we set  $\bcal_y=\{U(y,n):n<{\omega}\}$ then it follows that
$\bcal=\bigcup\{\bcal_y:y\in Y'\}$ is a 
 an  irreducible outer base of $Y$ in $X$ consisting of regular open sets. 
Now applying lemma \ref{lm:regirr} we can conclude the proof.
\eproof

Unfortunately, as $\cccmin \oo;$ is false, the above result is not applicable in the perhaps most interesting case when $\w(X)=\oo$. The annoying assumption
of  separability, however, can be circumvented as follows.

\begin{corollary}\label{cor:irred}
Assume $\cccmin {\kappa};$. 
If $X$ is a  first countable, 
regular space with $\w(X)\ge {\kappa}$, then there is an uncountable subspace 
$Y\subs X$ that has  an irreducible base. 
\end{corollary}

\proof
If $X$ is separable, then the previous theorem can be applied.
If $X$ is not separable, then $X$ contains an uncountable 
left separated subspace $Y$ and again $Y$ has an irreducible base.
\eproof

\end{document}